\documentclass[12pt]{article}

\usepackage{amsmath}
\allowdisplaybreaks[4]

\usepackage[all]{xy}
\usepackage{amssymb}
\usepackage{amsthm}
\usepackage{hyperref}
\hypersetup{colorlinks=true,linkcolor=blue,citecolor=red}
\usepackage{amsmath}
\usepackage{amscd}
\usepackage{verbatim}
\usepackage{eurosym}
\usepackage{float}
\usepackage{color}
\usepackage{dcolumn}
\usepackage[mathscr]{eucal}
\usepackage[all]{xy}
\usepackage{hyperref}
\usepackage{mathrsfs}
\usepackage{amsmath}
\usepackage{amssymb}
\usepackage{amsfonts,ifpdf}
\usepackage{graphicx}
\usepackage{times}
\usepackage{float}
\usepackage{epstopdf}
\usepackage{cite}
\usepackage{youngtab}
\usepackage{ytableau}
\ytableausetup
{mathmode, boxsize=0.9em}

\newlength{\boxedparwidth}
  {\begin{center} \begin{tabular}{|@{\hspace{.315in}}c@{\hspace{.15in}}|}
                  \hline \\ \begin{minipage}[t]{\boxedparwidth}
                  \setlength{\parindent}{.25in}}%
  {\end{minipage} \\ \\ \hline \end{tabular} \end{center}}

\allowdisplaybreaks

\newtheorem{thm}{Theorem}[section]
\newtheorem{prop}[thm]{Proposition}

\newtheorem{lem}[thm]{Lemma}
\newtheorem{cor}[thm]{Corollary}

\def\pf{\noindent{\it Proof.} }

\numberwithin{equation}{section}

\begin{document}

\begin{center}

 {\Large \bf Asymptotics for $k$-crank of $k$-colored partitions}
\end{center}
\vskip 3mm

\begin{center}
Helen W.J. Zhang$^{1,2}$ and Ying Zhong$^{1}$\\[8pt]
$^{1}$School of Mathematics\\
Hunan University\\
Changsha 410082, P. R. China\\[12pt]

$^{2}$Hunan Provincial Key Laboratory of \\
Intelligent Information Processing and Applied Mathematics\\
Changsha 410082, P. R. China\\[15pt]

Emails:  helenzhang@hnu.edu.cn,  YingZhong@hnu.edu.cn
\\[15pt]

\end{center}

\vskip 3mm

\begin{abstract}
In this paper, we obtain asymptotic formulas for $k$-crank of $k$-colored partitions.
Let $M_k(a, c; n)$ denote the number of $k$-colored partitions of $n$ with a $k$-crank congruent to $a$ mod $c$.
For the cases $k=2,3,4$, Fu and Tang derived several inequality relations for $M_k(a, c; n)$ using generating functions.
We employ the Hardy-Ramanujan Circle Method to extend the results of Fu and Tang.
Furthermore, additional inequality relations for $M_k(a, c; n)$ have been established, such as logarithmic concavity and logarithmic subadditivity.
\vskip 6pt

\noindent {\bf AMS Classifications}: 05A16, 05A17, 05A20, 11P82
\\ [7pt]
{\bf Keywords}: Asymptotics, Hardy-Ramanujan circle method, $k$-crank, $k$-colored partitions, Inequality

\end{abstract}

\section{Introduction}
The objective of this paper is to derive an asymptotic formula for the $k$-crank of $k$-colored partitions, and utilize the asymptotic formula to compute various inequality relations for the $k$-crank of $k$-colored partitions, such as unimodality, log-concavity, and log-subadditivity.
Recall that a partition of a positive integer $n$ is composed of a sequence of non-increasing positive integers that together sum to $n$. When each part in the partition can be represented by $k$ distinct colors, it is referred to as a $k$-colored partition.
Let $p(n)$ denote the number of partitions of $n$. Ramanujan provided three well-known congruences for $p(n)$.
To offer a combinatorial explanation for Ramanujan's congruences, the concepts of rank \cite{Atkin-Swinnerton-1954, Dyson-1944} and crank \cite{Andrews-Garvan-1988, Garvan-1988} were introduced within the framework of partition theory. In light of these developments, researchers have continued to explore various properties and relations associated with ranks and cranks. Let $N(a, c; n)$ denote the number of partitions of $n$ with rank congruent to $a$ modulo $c$. Bringmann \cite{Bringmann-2009} derived an asymptotic formula for $N(a, c; n)$ using the Hardy-Ramanujan circle method. Building on this research, Bringmann and Kane \cite{Bringmann-Kane-2010} subsequently provided an inequality for $N(a, c; n)$.
For further asymptotics and inequalities related to ranks and cranks, see, for example, \cite{Rolon-2015, Bringmann-Lovejoy-2007, Ciolan-2019, Ciolan-2022, Jennings-Shaffer-Reihill-2020}.

Building on earlier research in the domain of ranks and cranks, the exploration of $k$-colored partitions has emerged as a fresh area of interest. Motivated by Dyson's conjecture \cite{Dyson-1989}, the exploration of $k$-colored partitions has emerged as a fresh area of interest, leading to the introduction of the $k$-crank of $k$-colored partitions.
We employ the notation $C_k(x;q)$ from \cite{Bringmann-Dousse-2015} to signify the generating function of $M_k(m, n)$, as follows:
\begin{equation}\label{eq-ck}
C_k(x;q):=\sum_{n=0}^{\infty}\sum_{m=-\infty}^{\infty}M_k(m,n)x^mq^n=\frac{(q;q)_{\infty}^{2-k}}{(xq;q)_{\infty}(x^{-1}q;q)_{\infty}}.
\end{equation}
Fu and Tang \cite{Fu-Tang-2018} subsequently presented a combinatorial interpretation for $M_k(m, n)$, clarifying that $M_k(m, n)$ signifies the number of $k$-colored partitions of $n$ having a rank equal to $m$.
Acknowledging the importance of asymptotic properties, Bringmann and Dousse \cite{Bringmann-Dousse-2015} utilized the Wright Circle Method to unveil the asymptotic properties of $M_k(m, n)$.
Following this, by examining the generating function $C_k(x;q)$, Fu and Tang \cite{Fu-Tang-2018} managed to uncover a variety of fascinating properties and relationships for distinct values of $k$. The special case when $k=1$ corresponds to the crank of partitions. Moreover, it is worth noting that the particular instance of $k = 2$ results in the Hammond-Lewis birank \cite{Hammond-Lewis-2002} for $2$-colored partitions, while the case of $k = 3$ coincides with the multirank $r^*$ for $3$-colored partitions \cite{Fu-Tang-2018-2}. In this paper, we expand on these discoveries and employ the Hardy-Ramanujan Circle Method to determine the asymptotic properties of $M_k(a, c; n)$.

Before presenting our main result, we need to introduce some additional notations.
Let $0\leq h < p$ with $(h, p) = 1$, and define $h'$ such that $hh'\equiv-1~(\text{mod}~p)$ if $p$ is odd and $hh'\equiv-1~(\text{mod}~2p)$ if $p$ is even.
Furthermore, let $p_1:=\frac{p}{(p,c)},~c_1:=\frac{c}{(p,c)}$ and $0<l<c_1$ be defined by the congruence $l\equiv ap_1~(\text{mod}~c_1)$.
Lastly, consider $0 < a < c$ as coprime integers with $c$ being odd.
We now proceed to define the following sums for integers $m$ and $n$, both of which belong to $\mathbb{Z}$:

For $c\mid p$
\begin{align*}
 B_{a,c,p,k}(n,m):=(-1)^{ap+1}\sin\left(\frac{\pi a}{c}\right)\sum_{h~(\text{mod}~p)^*} \frac{\omega_{h,p}^k}{\sin\left(\frac{\pi ah'}{c}\right)} e^{-\frac{\pi ia^2p_1h'}{c}}e^{\frac{2\pi i}{p}\left(nh+mh'\right)},
\end{align*}
and for $c\nmid p$
\begin{align*}
 D_{a,c,p,k}(n,m):=(-1)^{ap+l}\sum_{h~(\text{mod}~p)^*}\omega_{h,p}^k e^{\frac{2\pi i}{p}\left(nh+m h'\right)}.
\end{align*}
In these expressions, the sums run through all primitive residue classes modulo $p$, and
\begin{equation*}
  \omega_{h,p}:=\exp\left(\pi i\sum_{\mu~~(\text{mod}~p)}\left(\left(\frac{\mu}{p}\right)\right)\left(\left(\frac{h\mu}{p}\right)\right)\right).
\end{equation*}
The sawtooth function is defined as follows:
\begin{equation*}
  \left(\left(x\right)\right):=
  \begin{cases}
  x-\lfloor x\rfloor-\frac 12, & \text{if}~ x\in\mathbb{R}\setminus \mathbb{Z},\\
  0, & \text{if}~x\in \mathbb{Z}.
  \end{cases}
\end{equation*}
Additionally, we have
\begin{align}
\delta_{a,c,p,r,k}^+&:=\frac{l^2}{2c_1^2}+\frac{k}{24}-\left(r+\frac 12\right)\frac{l}{c_1},
\label{delta+}\\[5pt]
\delta_{a,c,p,r,k}^-&:=\frac{l^2}{2c_1^2}+\frac{k}{24}-1-r\left(1-\frac l{c_1}\right)+\frac{l}{2c_1},
\label{delta-}
\end{align}
and
\begin{align*}
m_{a,c,p,r}^+&:=-\frac{a^2p_1p}{2cc_1}
+\frac{lap}{cc_1}-\frac{l^2}{2c_1^2}+\frac{rl}{c_1}-\frac{pra}{c}
+\frac{l}{2c_1}-\frac{ap}{2c},\\[5pt]
m_{a,c,p,r}^-&:=-\frac{a^2p_1p}{2cc_1}+\frac{lap}{cc_1}-\frac{l^2}{2c_1^2}
+r-\frac{rl}{c_1}+\frac{pra}{c}+1-\frac{l}{2c_1}+\frac{ap}{2c}.
\end{align*}

\begin{thm}\label{asym-M}
Let $\varepsilon > 0$ and consider $0 < a < c$ as coprime integers. Given that $c$ is odd and $n$ is a positive integer, the following holds for $k \leq 12$:
\begin{align}\label{eq-asym-M}
M_k(a,c;n)=
&\frac{ 2\pi}{c} \left(n-\frac{k}{24}\right)^{-\frac k4-\frac 12}\left(\frac{k}{24}\right)^{\frac k4+\frac 12}\sum_{p=1}^{\infty}\frac{A_{p,k}(n,0)}{p}
I_{\frac k2+1}\left(\frac{4\pi}{p}\sqrt{\frac{k}{24}\left(n-\frac{k}{24}\right)}\right)
\notag\\[5pt]
&+\frac{1}{c}\sum_{\beta=1}^{c-1}\zeta_c^{-a\beta}
\left[2\pi i\left(\frac{k}{24n-k}\right)^{\frac k4}\sum_{1\leq p\leq N\atop c\mid p}\frac{B_{\beta,c,p,k}(-n,0)}{p}I_{\frac k2}\left(\frac{\pi\sqrt{k(24n-k)}}{6p}\right)\right.
\notag\\[5pt]
&+\frac{4 \cdot24^{\frac k4}\pi\sin\left(\frac{\pi \beta}{c}\right)}{\left(24n-k\right)^{\frac k4}}
  \sum_{\substack{p,r\\ c\nmid p\\ \delta_{\beta,c,p,r,k}^j>0 \\j\in\{+,-\} }} \frac{\left(\delta_{\beta,c,p,r,k}^j\right)^{\frac k4}D_{\beta,c,p,k}(-n,m_{\beta,c,p,r}^j)}{p}
\notag\\[5pt]
&\left.\times I_{\frac k2}\left(\frac{\pi}{p}\sqrt{\frac{2 \delta_{\beta,c,p,r,k}^j(24n-k)}{3}}\right)\right]+O(n^\varepsilon),
\end{align}
where $p_k(n)$ is the number of $k$-colored partitions of $n$.
\end{thm}

We find the main term of $M_k(a,c;n)$ in Theorem \ref{asym-M}, which leads us to the following asymptotic formula.
\begin{thm}\label{c-asym}
Let $0 < a < c$ be coprime integers. Given that $c$ is odd and $n$ is a positive integer, the following relation holds for $k \leq 12$:
\begin{align*}
M_k(a,c;n) \sim \frac{p_k(n)}{c} \sim \frac{2}{c}\left(\frac{k}{3}\right)^{\frac{1+k}{4}}(8n)^{-\frac{3+k}{4}}e^{\pi\sqrt{\frac{2kn}{3}}},
\end{align*}
as $n \rightarrow \infty$.
\end{thm}

Asymptotic formulas for rank or crank often give rise to rank inequalities, as demonstrated in the works of Bringmann and Kane \cite{Bringmann-Kane-2010}, Hamakiotes, Kriegman and Tsai \cite{Hamakiotes-Kriegman-Tsai-2021}, Zhang and Zhong \cite{Zhang-Zhong-2023}.
In this paper, we extend the inequalities of Fu and Tang \cite{Fu-Tang-2018} for $M_k(a, c; n)$, as presented in Theorem \ref{ineq-M3}. Furthermore, we establish the strict log-subadditivity and log-concavity of $M_k(a, c; n)$ using the above asymptotic formula.
\begin{thm}\label{ineq-M3}
Let $0\leq a<b\leq \frac{c-1}{2}$ . Given that $c\geq 2$ (resp. $c\geq 4$) and $5\leq k\leq 12$ (resp. $k=3$) are odd integers, then for $n > N_{a,b,c,k}$, we have
\begin{equation*}
  M_k(a,c;n)>M_k(b,c;n),
\end{equation*}
where $N_{a,b,c,k}$ is an explicit constant.
\end{thm}

\begin{thm}\label{addi}
Let $0<a<c$ be coprime integers. Given that $c$ is odd, then for $k\leq12$ we have
\begin{equation*}
  M_k(a,c;n_1)M_k(a,c;n_2)>M_k(a,c;n_1+n_2)
\end{equation*}
for sufficiently large $n_1$ and $n_2$.
\end{thm}

\begin{thm}\label{ineq-concave}
Let $0<a<c$ be coprime integers. Given that $c$ is odd, then for $k\leq12$ we have
\begin{equation*}
  M_k(a,c;n_1)M_k(a,c;n_2)>M_k(a,c;n_1-1)M_k(a,c;n_2+1)
\end{equation*}
for sufficiently large $n_1$ and $n_2$ with $n_1<n_2+1$.
\end{thm}

By setting $n_1 = n_2$ in Theorem \ref{ineq-concave}, we can deduce that $M_k(a,c;n)$ is asymptotically log-concave.
\begin{cor}
Let $0<a<c$ be coprime integers. Given that $c$ is odd and $n$ is a positive integer, then for $k\leq12$ and sufficiently large $n$ we have
\begin{equation*}
  M_k(a,c;n)^2>M_k(a,c;n-1)M_k(a,c;n+1).
\end{equation*}
\end{cor}

In Theorem \ref{addi}, we demonstrate the strict log-subadditivity of $M_k(a,c;n)$ for sufficiently large $n_1$ and $n_2$. We also provide an example to clarify the exact bounds of $n_1$ and $n_2$ in relation to Theorem \ref{asym-M}.
\begin{thm}\label{ineq-sum}
Let $c$ be an odd integer and define
$$M_c:=\frac{2^3\times3^{12}\times(4.0102\times 10^6c^{\frac 52}+1824.9112c)^4}{\left(\pi-2\pi\sqrt{\frac{1}{c^2}-\frac{1}{c}+\frac{1}{4}}\right)^{12}}+1,$$
then we have
\begin{equation}\label{thm-ineq}
  M_3(a,c;n_1)M_3(a,c;n_2)>M_3(a,c;n_1+n_2)
\end{equation}
for $c\geq 4$ and $n_1,n_2\geq M_c$ or $c=3$ and $n_1,n_2\geq 4$.
\end{thm}

The organization of this paper is as follows. In Section \ref{trans-sec}, we present a transformation formula for $C_k(x;q)$. In Section \ref{thm}, we prove a proposition using the Hardy-Ramanujan circle method, which enables us to further prove Theorems \ref{asym-M} and \ref{c-asym}.
Section \ref{ineq} presents proofs of several inequality relations for $M_k(a,c;n)$, including unimodality, log-concavity, and log-subadditivity.
In Section \ref{SLS-3}, we derive an upper bound and a lower bound for $M_3(a,c;n)$. With their assistance, we demonstrate the strict log-subadditivity of $M_3(a,c;n)$.

\section{A transformation formula}\label{trans-sec}
In this section, our primary objective is to establish a transformation formula for the function $C_k(x;q)$. To achieve this goal, we first introduce some notations and auxiliary functions.

Take into account the following representation for the complex variables $\tau_1$ and $\tau_2$:
\begin{align*}
\tau_1:=\frac{h+iz}{p},\qquad \tau_2:=\frac{h'+\frac{i}{z}}{p},
\end{align*}
and similarly for $\omega_1$, $\omega_2$, and $\omega_3$:
\begin{align*}
\omega_1:=\frac{a}{c},\qquad \omega_2:=\frac{ah'}{c},\qquad \omega_3:=-h'\omega_1+\frac{l}{c_1}\tau_2.
\end{align*}
We define the function $\chi(h, h', p)$ as follows:
\begin{equation*}
\chi(h, h', p):=i^{-\frac 12}\omega_{h,p}^{-1}e^{-\frac {\pi i}{12p}(h'-h)}.
\end{equation*}
Here, $h'$ is a solution to the congruence relation $hh' \equiv -1 \pmod{p}$.

The Dedekind $\eta$-function is given by the expression
\begin{equation}\label{eta-d}
\eta(\tau):=q^{\frac 1{24}}\left(q;q\right)_\infty,
\end{equation}
where $q:=e^{2\pi i\tau}$ and $\tau\in\mathbb{H}$. For $z\in \mathbb{C}$ with $\text{Re}(z)>0$, the $\eta$-function satisfies the relation
\begin{equation}\label{eta-t}
\eta\left(\tau_1\right)=\sqrt{\frac iz}\chi(h,h',p) \eta\left(\tau_2\right),
\end{equation}
where we take the principal branch of the square root.

Additionally, the Jacobi $\vartheta$-function is defined as
\begin{equation*}
\vartheta(u;\tau):=\sum_{\nu \in\mathbb{Z}+\frac 12}e^{\pi i\nu^2\tau+2\pi i\nu\left(u+\frac 12\right)}.
\end{equation*}
We define $x:=e^{2\pi i\omega}$, where $\omega\in \mathbb{C}$. Let $h, p$ be coprime integers, with $h'$ defined as mentioned earlier. The $\vartheta$-function satisfies the following properties:
\begin{equation}\label{theta-t1}
\vartheta(\omega;\tau)=-2\sin(\pi\omega)q^{\frac 18}(q)\infty(xq)\infty(x^{-1}q)_\infty,
\end{equation}
and for $\text{Re}(z)>0$,
\begin{equation}\label{theta-t2}
\vartheta\left(\omega;\tau_1\right)=\chi(h,h',p)^3\sqrt{\frac iz}e^{-\frac{\pi p\omega^2}{z}} \vartheta\left(\frac{i\omega}{z};\tau_2\right).
\end{equation}
According to the works of \cite{Andrews-Garvan-1988, Dyson-1944, Garvan-1988}, we know that
\begin{equation*}
\frac{(q;q)_{\infty}}{(xq;q)_{\infty}(x^{-1}q;q)_{\infty}}
=\frac{1-x}{(q;q)_{\infty}}\sum_{n\in\mathbb{Z}}\frac{(-1)^nq^{\frac {n(n+1)}2}}{1-xq^n}.
\end{equation*}
Multiplying both sides of the above equation by $(q;q)_{\infty}^{1-k}$, combining the definition of $C_k(x;q)$ \eqref{eq-ck}, we have
\begin{equation}\label{C-k}
C_k(x;q)
=\frac{1-x}{(q;q)_{\infty}^k}\sum_{n\in\mathbb{Z}}\frac{(-1)^nq^{\frac {n(n+1)}2}}{1-xq^n}.
\end{equation}
By combining equations \eqref{eta-d}, \eqref{theta-t1} and \eqref{C-k}, we arrive at the following expression for $C_k(x;q)$:
\begin{equation}\label{C-k-t}
C_k(x;q)=\frac {-2\sin(\pi\omega)q^{\frac {k}{24}}\eta(\tau)^{3-k}}{\vartheta(\omega;\tau)}.
\end{equation}
With this result in hand, we now proceed to construct the transformation formula for $C_k(x;q)$ using the introduced notations and auxiliary functions.

\begin{prop}\label{trans}
Let $q_1:=e^{2\pi i\tau_2}$, we have\\
(1)~for $c\mid p$, then
\begin{align*}
  C_k\left(e^{2\pi i\omega_1};e^{2\pi i\tau_1}\right)=&\frac{i\sin\left(\pi\omega_1\right)z^{\frac k2-1}}{\sin\left(\pi\omega_2\right)}(-1)^{ap+1}\omega_{h,p}^k e^{\frac{k\pi}{12p}(z^{-1}-z)-\frac{\pi ia^2p_1h'}{cc_1}}\\
  &\times C_k\left(e^{2\pi i\omega_2};e^{2\pi i\tau_2}\right),
\end{align*}
(2)~for $c\nmid p$, then
\begin{align*}
  C_k\left(e^{2\pi i\omega_1};e^{2\pi i\tau_1}\right)=&4i(-1)^{ap+l+1}\sin\left(\pi\omega_1\right)z^{\frac k2-1}\omega_{h,p}^k
  e^{-\frac{\pi ia^2p_1h'}{cc_1}+\frac{2\pi ialh'}{cc_1}}\\
  &\times e^{\frac{k\pi}{12p}(z^{-1}-z)}q_1^{-\frac{l^2}{2c_1^2}} C_k\left(ah',\frac {lc}{c_1},c;q_1\right),
\end{align*}
where
\begin{equation}\label{C-k-4v}
  C_k(a,b,c;q):=\frac{i}{2(q)_{\infty}^k}\sum_{m=-\infty}^{\infty}\frac{(-1)^m e^{-\frac{\pi ia}{c}}q^{\frac{m(m+1)}{2}+\frac b{2c}}}{1-e^{-\frac{2\pi ia}{c}}q^{m+\frac bc}}.
\end{equation}
\end{prop}
\pf
(1). By substituting $\tau=\tau_1$ and $\omega=\omega_1$ into equation \eqref{C-k-t} and combining equations \eqref{eta-t} and \eqref{theta-t2}, we obtain:
\begin{equation}\label{c-k}
C_k\left(e^{2\pi i\omega_1};e^{2\pi i\tau_1}\right)=-2\sin(\pi \omega_1)iz^{\frac k2-1}\omega_{h,p}^k e^{\frac{k\pi i}{12p}(h'+iz)+\frac{\pi p\omega_1^2}{z}}\frac{\eta\left(\tau_2\right)^{3-k}}{\vartheta\left(\frac{i\omega_1}{z};\tau_2\right)}.
\end{equation}
Next, set $\tau=\tau_2$ and $\omega=\omega_2$ in equation \eqref{C-k-t} to arrive at:
\begin{align}\label{eta^3}
\eta\left(\tau_2\right)^{3-k}=\frac{ C_k\left(e^{2\pi i\omega_2};e^{2\pi i\tau_2}\right)\vartheta\left(\omega_2;\tau_2\right)}{-2\sin\left(\pi\omega_2\right)e^{\frac{k\pi i\tau_2}{12}}}.
\end{align}
Combining this with equation \eqref{eta^3}, we then obtain:
\begin{align}\label{c-k-1}
C_k\left(e^{2\pi i\omega_1};e^{2\pi i\tau_1}\right)=&\frac{i\sin\left(\pi\omega_1\right)z^{\frac k2-1}}{\sin\left(\pi\omega_2\right)}\omega_{h,p}^k e^{\frac{k\pi}{12p}(z^{-1}-z)+\frac{\pi a^2p}{zc^2}}
\nonumber\\[5pt]
&\times\frac{\vartheta\left(\omega_2;\tau_2\right)}
{\vartheta\left(\frac{i\omega_1}{z};\tau_2\right)}
 C_k\left(e^{2\pi i\omega_2};e^{2\pi i\tau_2}\right).
\end{align}
Since $c\mid p$, $A:=\frac{ap}{c}\in\mathbb{Z}$, Rolon \cite[p.153]{Rolon-2015} obtained
\begin{equation*}
\frac{\vartheta\left(\omega_2;\tau_2\right)}
{\vartheta\left(\frac{i\omega_1}{z};\tau_2\right)}=(-1)^{ap+1}e^{-\frac{\pi ia^2ph'}{c^2}}e^{-\frac{\pi a^2p}{zc^2}}.
\end{equation*}
Inserting this expression to \eqref{c-k-1} yields the transformation formula for the case of $c\mid p$.

(2). For $c \nmid p$, $\frac{ap}{c} \notin \mathbb{Z}$. We need to explore another approach to obtain the desired transformation formula. Here, we define $l$ as a solution to the congruence condition $l \equiv ap_1 \pmod{c_1}$. Consequently, it is evident that $B := \frac{l - ap_1}{c_1} \in \mathbb{Z}$. In this context, Rolon \cite[p.153]{Rolon-2015} provided the following equation:
\begin{equation*}
\vartheta\left(\frac{i\omega_1}{z};\tau_2\right)=(-1)^{ap+l}e^{\frac{\pi i(l-ap_1)^2}{c_1^2}\tau_2}e^{ \frac{2\pi i(l -ap_1)}{c_1}\cdot\frac{i\omega_1}{z}}
\vartheta\left(\omega_3;\tau_2\right).
\end{equation*}
Together with \eqref{c-k}, we have
\begin{align*}
  C_k\left(e^{2\pi i\omega_1};e^{2\pi i\tau_1}\right)=&2i\sin\left(\pi\omega_1\right)z^{\frac k2-1}(-1)^{ap+l+1}\omega_{h,p}^k e^{-\frac{\pi ia^2p_1h'}{cc_1}+\frac{2\pi ialh'}{cc_1}+\frac{k\pi i}{12p}(h'+iz)} q_1^{-\frac{l^2}{2c_1^2}}
  \frac{\eta\left(\tau_2\right)^{3-k}}{\vartheta\left(\omega_3;\tau_2\right)}.
\end{align*}
Now setting $\tau=\tau_2$ and $\omega=\omega_3$ in \eqref{C-k-t}, the above equation can be rewritten as
\begin{align}\label{c-k-2}
  C_k\left(e^{2\pi i\omega_1};e^{2\pi i\tau_1}\right)=&i\sin\left(\pi\omega_1\right)z^{\frac k2-1}(-1)^{ap+l}\omega_{h,p}^k
\nonumber\\[5pt]
&\times e^{\frac{k\pi}{12p}(z^{-1}-z)-\frac{\pi ia^2p_1h'}{cc_1}+\frac{2\pi ialh'}{cc_1}} q_1^{-\frac{l^2}{2c_1^2}}\frac{C_k\left(e^{2\pi i\omega_3};q_1\right)}{\sin\left(\pi\omega_3\right)}.
\end{align}
Utilizing the exponential representation of the sine function and \eqref{C-k}, we deduce:
\begin{align*}
\frac{C_k\left(e^{2\pi i\omega_3};q_1\right)}{\sin\left(\pi\omega_3\right)}&=
 \frac{2i \left(1-e^{2\pi i\omega_3}\right)}{(q_1)_{\infty}^k\left(e^{\pi i\omega_3}-e^{-\pi i\omega_3}\right)}\sum_{m\in\mathbb{Z}}\frac{(-1)^mq_1^{m(m+1)/2}}{1-e^{2\pi i\omega_3}q_1^m}\\
 &=\frac{-2i e^{\pi i\omega_3}}{(q_1)_{\infty}^k}\sum_{m\in\mathbb{Z}}\frac{(-1)^mq_1^{m(m+1)/2}}{1-e^{2\pi i\omega_3}q_1^m}.
\end{align*}
By substituting the above identity into \eqref{c-k-2}, we obtain:
\begin{align*}
C_k\left(e^{2\pi i\omega_1};e^{2\pi i\tau_1}\right)
  &=i\sin\left(\pi\omega_1\right)z^{\frac k2-1}(-1)^{ap+l+1}\omega_{h,p}^k e^{\frac{k\pi}{12p}(z^{-1}-z)-\frac{\pi ia^2p_1h'}{cc_1}+\frac{2\pi ialh'}{cc_1}} q_1^{-\frac{l^2}{2c_1^2}}
  \\
  &\quad\times \frac{2i e^{\pi i\omega_3}}{(q_1)_{\infty}^k}\sum_{m\in\mathbb{Z}}\frac{(-1)^mq_1^{m(m+1)/2}}{1-e^{2\pi i\omega_3}q_1^m}\\
  &=i\sin\left(\pi\omega_1\right)z^{\frac k2-1}(-1)^{ap+l+1}\omega_{h,p}^k e^{\frac{k\pi}{12p}(z^{-1}-z)-\frac{\pi ia^2p_1h'}{cc_1}+\frac{2\pi ialh'}{cc_1}} q_1^{-\frac{l^2}{2c_1^2}}
  \\
  &\quad\times\frac{2i }{(q_1)_{\infty}^k}\sum_{m\in\mathbb{Z}}\frac{(-1)^m e^{-\pi ih'\omega_1}q_1^{\frac l{2c_1}+\frac{m(m+1)}2}}{1-e^{-2\pi ih'\omega_1}q_1^{m+\frac l{c_1}}},
\end{align*}
which completes the proof.
\qed

\section{Asymptotic formula for $M_k(a,c;n)$}\label{thm}
The following theorem plays a crucial role in proving the asymptotic formula.
\begin{lem}\label{bound}
Let $n, m, p, D \in \mathbb{Z}$ with $(D, p) = 1$, and $0 \leq \sigma_1 < \sigma_2 \leq p$. Then there exist constants $C_1$ and $C_2$ such that:

\rm{(1)} We have
\begin{align}\label{asym-form-1}
|F|:=\left|\sum_{h\pmod{p}^*\atop \sigma_1\leq Dh'\leq\sigma_2}
\omega_{h,p}^k e^{\frac{2\pi i}{p}(-hn+h'm)}
\right|\leq C_1\cdot(24n+k,p)^{\frac{1}{2}}p^{\frac{1}{2}+\varepsilon}.
\end{align}

\rm{(2)} We have
\begin{align}\label{asym-form-2}
\left|\frac{\sin\left(\frac{\pi a}{c}\right)}{(-1)^{ap+1}}
\sum_{h\pmod{p}^*\atop \sigma_1\leq Dh'\leq\sigma_2}\frac{\omega_{h,p}^k}{\sin\left(\frac{\pi ah'}{c}\right)}e^{-\frac{\pi ia^2p_1h'}{c}}e^{\frac{2\pi i}{p}(-hn+h'm)}
\right|\leq C_2\cdot(24n+k,p)^{\frac{1}{2}}p^{\frac{1}{2}+\varepsilon}.
\end{align}
The constants $C_1$ and $C_2$ are independent of $a$ and $p$.
\end{lem}
\pf
According to the calculations of Andrews \cite[p. 482]{Andrews-1966} and defining $\theta := (h, 24)$, we have
\begin{align*}
  \omega_{h,p}^k=(-i)^k\left(\frac{-p}{h} \right)^k\exp\left( 2k\pi i\frac{(4p^2+3p-1)h-(p^2-1)h'}{24p}\right),
\end{align*}
when $p$ is even and $hh'\equiv-1\pmod{24p/\theta}$.
Additionally,
\begin{align*}
  \omega_{h,p}^k=\exp\left( \frac{-k\pi i(p-1)}{4}\right)\left(\frac{-h}{p} \right)^k\exp\left( 2k\pi i\frac{(p^2-1)h-(p^2-1)h'}{24p}\right)
\end{align*}
when $p$ is odd  and $hh'\equiv-1\pmod{24p/\theta}$.

Hence, $F$ can be expressed as:
\begin{align*}
 \sum_{h\pmod{p}^*\atop \sigma_1\leq Dh'\leq\sigma_2}(-i)^k\left(\frac{-p}{h} \right)^k\exp\left( 2\pi i\frac{\left((4p^2+3p-1)k-24n\right)h-\left((p^2-1)k-24m\right)h'}{24p}\right)
\end{align*}
when $p$ is even and $hh'\equiv-1\pmod{24p/\theta}$.
For odd $p$ and $hh' \equiv -1 \pmod{24p/\theta}$, it can be written as:
\begin{align*}
 \sum_{h\pmod{p}^*\atop \sigma_1\leq Dh'\leq\sigma_2}\exp\left( \frac{-k\pi i(p-1)}{4}\right)&\left(\frac{-h}{p} \right)^k
 \\
 &\times\exp\left( 2\pi i\frac{\left((p^2-1)k-24n\right)h-\left((p^2-1)k-24m\right)h'}{24p}\right).
\end{align*}
We have
\begin{align*}
  |F|=O\left(\sum_{h\pmod{p}^*\atop \sigma_1\leq Dh'\leq\sigma_2}\exp\left( 2\pi i\frac{\left((4p^2+3p-1)k-24n\right)h-\left((p^2-1)k-24m\right)h'}{24p}\right) \right)
\end{align*}
for even $p$.
When $p$ is odd,
\begin{align*}
  |F|=O\left( \sum_{h\pmod{p}^*\atop \sigma_1\leq Dh'\leq\sigma_2}\exp\left( 2\pi i\frac{\left((p^2-1)k-24n\right)h-\left((p^2-1)k-24m\right)h'}{24p}\right)\right)
\end{align*}
By utilizing an estimate of Sali'{e} \cite[eq. (5)]{Salie-1933}, we have
\begin{align*}
 |F|&=\begin{cases}
 O\left(p^{\frac 12+\varepsilon}\left( (4p^2+3p-1)k-24n,24p\right) \right),& p~ \text{even},\\[5pt]
 O\left(p^{\frac 12+\varepsilon}\left( (p^2-1)k-24n,24p\right) \right), & p~ \text{odd},
 \end{cases}\\[5pt]
 &=\begin{cases}
 O\left(p^{\frac 12+\varepsilon}\left( 24n+k,p\right) \right),& p~ \text{even},\\[5pt]
 O\left(p^{\frac 12+\varepsilon}\left( 24n+k,p\right) \right), & p~ \text{odd}.
 \end{cases}
\end{align*}
We have verified that \eqref{asym-form-2} holds. We define $\overline{c} := c$ if $p$ is odd and $\overline{c} := 2c$ if $p$ is even. It is clear that $\frac{e^{-\frac{\pi ia^2p_1h'}{c}}}{\sin\left(\frac{\pi ah'}{c}\right)}$ depends only on the residue class of $h'(\text{mod}~\overline{c})$. We can then express \eqref{asym-form-2} as:
\begin{equation*}
  \frac{\sin\left(\frac{\pi a}{c}\right)}{(-1)^{ap+1}}\sum_{c_j}\frac{e^{-\frac{\pi ia^2p_1h'}{c}}}{\sin\left(\frac{\pi ah'}{c}\right)}
\sum_{\substack{h\pmod{p}^*\\ \sigma_1\leq Dh'\leq\sigma_2\\ h'\equiv c_j~(\text{mod}~ \overline{c} )}}\omega_{h,p}^ke^{\frac{2\pi i}{p}(-hn+h'm)},
\end{equation*}
where $c_j$ runs through a set of primitive residues $(\text{mod}~\overline{c})$. The inner sum can be rewritten as
\begin{align*}
  \frac 1{\overline{c}}\sum_{\substack{h\pmod{p}^*\\ \sigma_1\leq Dh'\leq\sigma_2}}\omega_{h,p}^ke^{\frac{2\pi i}{p}(-hn+h'm)}
  &\sum_{r~~(\text{mod}~ \overline{c})}e^{\frac{2\pi ir}{\overline{c}}(h'-c_j)}\\
  &=\frac 1{\overline{c}}\sum_{r~~(\text{mod}~ \overline{c})}e^{-\frac{2\pi irc_j}{\overline{c}}}\sum_{\substack{h\pmod{p}^*\\ \sigma_1\leq Dh'\leq\sigma_2}}\omega_{h,p}^ke^{\frac{2\pi i}{p}(-hn+h'(m+\frac{pr}{\overline{c}}))}.
\end{align*}
Applying \eqref{asym-form-1} effortlessly yields \eqref{asym-form-2}, thus completing the proof.
\qed

Consider the expression
\[C_k\left(e^{\frac{2\pi ia}{c}};q\right)
=1+\sum_{n=1}^\infty A_k\left(\frac{a}{c};n\right)q^n.\]
We derive the subsequent asymptotic formulas for the coefficient $A_k\left(\frac{a}{c};n\right)$.

\begin{prop}\label{asym}
Given $\varepsilon > 0$, suppose $0 < a < c$ are coprime integers, $c$ is odd, and $n$ is a positive integer. Then, for $k \leq 12$, we obtain
\begin{align*}
A_k\left(\frac{a}{c};n\right)
&=2\pi i\left(\frac{k}{24n-k}\right)^{\frac k4}\sum_{1\leq p\leq N\atop c\mid p}\frac{B_{a,c,p,k}(-n,0)}{p}I_{\frac k2}\left(\frac{\pi\sqrt{k(24n-k)}}{6p}\right)+\frac{4 \cdot24^{\frac k4}\pi\sin\left(\frac{\pi a}{c}\right)}{\left(24n-k\right)^{\frac k4}}
\\[5pt]
&\quad\times \sum_{\substack{p,r\\ c\nmid p\\ \delta_{a,c,p,r,k}^j>0 \\j\in\{+,-\} }} \frac{\left(\delta_{a,c,p,r,k}^j\right)^{\frac k4}D_{a,c,p,k}(-n,m_{a,c,p,r}^j)}{p}I_{\frac k2}\left(\frac{\pi}{p}\sqrt{\frac{2 \delta_{a,c,p,r,k}^j(24n-k)}{3}}\right)+O(n^\varepsilon).
\end{align*}
\end{prop}

\pf
Applying Cauchy's theorem, we obtain for $n > 0$
\[A_k\left(\frac{a}{c};n\right)=\frac{1}{2\pi i}\int_\mathcal{C}\frac{C_k\left(e^{\frac{2\pi ia}{c}};q\right)}{q^{n+1}}dq,\]
where $\mathcal{C}$ represents an arbitrary path inside the unit circle encircling $0$ counterclockwise. Choosing a circle with radius $e^{-\frac{2\pi}{n}}$ and parametrization $q=e^{-\frac{2\pi}{n}+2\pi it}$ for $0\leq t\leq1$, we get
\[A_k\left(\frac{a}{c};n\right)=\int_0^1C_k\left(e^{\frac{2\pi ia}{c}};e^{-\frac{2\pi}{n}+2\pi it}\right)e^{2\pi-2\pi int}dt.\]
We define
\[\vartheta'_{h,p}:=\frac{1}{p\left(\widetilde{p}_1+p\right)},~~~~~
\vartheta''_{h,p}:=\frac{1}{p\left(\widetilde{p}_2+p\right)},\]
where $\frac{h_1}{\widetilde{p}_1}<\frac{h}{p}<\frac{h_2}{\widetilde{p}_2}$ are adjacent Farey fractions in the Farey sequence of order $N:=\lfloor n^{1/2}\rfloor$.
From the theory of Farey fractions, we know that
\[\frac{1}{p+\widetilde{p}_j}\leq\frac{1}{N+1}~~~~(j=1,2).\]
Now, we decompose the path of integration along Farey arcs $-\vartheta'_{h,p}\leq\Phi\leq\vartheta''_{h,p}$, where $\Phi=t-\frac{h}{p}$ and $0\leq h< p\leq N$ with $(h,p)=1$. From this decomposition of the path we can rewrite the integral along these arcs:
\[A_k\left(\frac{a}{c};n\right)=\sum_{h,p}e^{-\frac{2\pi ihn}{p}}
\int_{-\vartheta'_{h,p}}^{\vartheta''_{h,p}}C_k\left(e^{\frac{2\pi ia}{c}};e^{\frac{2\pi i}{p}(h+iz)}\right)e^{\frac{2\pi nz}{p}}d\Phi,\]
where $z=\frac{p}{n}-p\Phi i$.
We insert our transformation formula from Proposition \ref{trans} into the integral and obtain
\begin{align*}
A_k\left(\frac{a}{c};n\right)
&=i\sin\left(\frac{\pi a}{c}\right)\sum_{h,p\atop c\mid p}\omega_{h,p}^k\frac{(-1)^{ap+1}}{\sin\left(\frac{\pi ah'}{c}\right)} e^{-\frac{\pi ia^2p_1h'}{c}-\frac{2\pi ihn}{p}}
\\
&\quad\times\int_{-\vartheta'_{h,p}}^{\vartheta''_{h,p}} z^{\frac k2-1}e^{\frac {2\pi z}{p}\left(n-\frac {k}{24}\right)+\frac{k\pi}{12pz}}
C_k\left(e^{\frac{2\pi iah'}{c}};q_1\right)d\Phi
\\
&\quad -4i\sin\left(\frac{\pi a}{c}\right)\sum_{h,p\atop c\nmid p}\omega_{h,p}^k(-1)^{ap+l} e^{-\frac{\pi ia^2p_1h'}{cc_1}+\frac{2\pi ialh'}{cc_1}-\frac{2\pi ihn}{p}}
\\
&\quad\times \int_{-\vartheta'_{h,p}}^{\vartheta''_{h,p}} z^{\frac k2-1}e^{\frac {2\pi z}{p}\left(n-\frac {k}{24}\right)+\frac{k\pi}{12pz}} q_1^{-\frac{l^2}{2c_1^2}} C_k\left(ah',\frac {lc}{c_1},c;q_1\right),
\end{align*}
and we call them $\Sigma_1$ and $\Sigma_2$ respectively.
To estimate $\Sigma_1$, we first consider the principal part of $C_k\left(e^{\frac{2\pi iah'}{c}};q_1\right)$ with respect to the $q_1$ variable. Using \eqref{eq-ck}, we obtain
\begin{align*}
 C_k\left(e^{\frac{2\pi iah'}{c}};q_1\right)&=1+\sum_{r=1}^\infty\sum_{m=-\infty}^{\infty}M_k(m,r)e^{\frac{2\pi iah'm}c}q_1^r\\
 &=1+\sum_{r=1}^\infty\sum_{s~~(\text{mod}~c)}\sum_{t=-\infty}^{\infty}M_k(s+tc,r)e^{\frac{2\pi iah'm}c}q_1^r.
\end{align*}
Letting \[a_k(r,s):=\sum_{t=-\infty}^{\infty}M_k(s+tc,r),\] we obtain
\begin{align}\label{C-k-2}
C_k\left(e^{\frac{2\pi iah'}{c}};q_1\right)=1+\sum_{r=1}^\infty\sum_{s~~(\text{mod}~c)}a_k(r,s)e^{\frac{2\pi ih'}cm_{r,s}}q_1^r,
\end{align}
where $m_{r,s}\in\mathbb{Z}$ and \[\sum_{s~(\text{mod}~c)}a_k(r,s)=p_k(r)\] for $r\geq1$.
Therefore, $\Sigma_1$ part can be expressed as
\[\Sigma_1=S_1+S_2.\]
Here
\begin{align}\label{S1}
  S_1=i\sin\left(\frac{\pi a}{c}\right)\sum_{h,p\atop c\mid p}\omega_{h,p}^k\frac{(-1)^{ap+1}}{\sin\left(\frac{\pi ah'}{c}\right)} e^{-\frac{\pi ia^2p_1h'}{c}-\frac{2\pi ihn}{p}}\int_{-\vartheta'_{h,p}}^{\vartheta''_{h,p}} z^{\frac k2-1}e^{\frac {2\pi z}{p}\left(n-\frac {k}{24}\right)+\frac{k\pi}{12pz}}d\Phi
\end{align}
and
\begin{align}\label{S2}
  S_2=&i\sin\left(\frac{\pi a}{c}\right)\sum_{h,p\atop c\mid p}\omega_{h,p}^k\frac{(-1)^{ap+1}}{\sin\left(\frac{\pi ah'}{c}\right)} e^{-\frac{\pi ia^2p_1h'}{c}-\frac{2\pi ihn}{p}}
  \notag\\
  &\times\int_{-\vartheta'_{h,p}}^{\vartheta''_{h,p}} z^{\frac k2-1}e^{\frac {2\pi z}{p}\left(n-\frac {k}{24}\right)+\frac{k\pi}{12pz}}\sum_{r=1}^\infty p_k(r)e^{\frac{2\pi ih'}cm_{r,s}}q_1^rd\Phi.
\end{align}
In order to bound $S_2$, we split the integral in the following way
\[\int_{-\vartheta'_{h,p}}^{\vartheta''_{h,p}}
=\int_{-\frac{1}{p(N+p)}}^{\frac{1}{p(N+p)}}+\int_{-\frac{1}{p(\widetilde{p}_1+p)}}^{-\frac{1}{p(N+p)}}
+\int_{\frac{1}{p(N+p)}}^{\frac{1}{p(\widetilde{p}_2+p)}},\]
Consequently, $S_2$ can be decomposed into three integrals:
\[S_2=S_{21}+S_{22}+S_{23}.\]
For example,
\begin{align*}
  S_{21}=&i\sin\left(\frac{\pi a}{c}\right)\sum_{h,p\atop c\mid p}\omega_{h,p}^k\frac{(-1)^{ap+1}}{\sin\left(\frac{\pi ah'}{c}\right)} e^{-\frac{\pi ia^2p_1h'}{c}-\frac{2\pi ihn}{p}}\\
  &\times\int_{-\frac{1}{p(N+p)}}^{\frac{1}{p(N+p)}} z^{\frac k2-1}e^{\frac {2\pi z}{p}\left(n-\frac {k}{24}\right)+\frac{k\pi}{12pz}}\sum_{r=1}^\infty p_k(r)e^{\frac{2\pi ih'}cm_{r,s}}q_1^rd\Phi.
\end{align*}
One can easily verify that $\frac pn\leq|z|\leq\sqrt{\frac 2n}$, hence
\begin{equation*}
|z|^{\frac k2-1}\leq
\begin{cases}
p^{-\frac 12}n^{\frac 12},&k=1,\\
\left(\frac 2n\right)^{\frac k4-\frac 12},&2\leq k\leq12.
\end{cases}
\end{equation*}
Taking the absolute value of $S_{21}$, it is same as the result of Rolon \cite[p.158]{Rolon-2015} for $k=1$, i.e., $|S_{21}|=O(n^\varepsilon)$. By \cite[p.492]{Dragonette-1952}, we have $\operatorname{Re}(\frac 1z)<p$. When $2\leq k\leq12$, using the fact that $\operatorname{Re}(z)=\frac pn$ and $\frac p2\leq \operatorname{Re}(\frac 1z)<p$ then
\begin{align*}
 \left|S_{21}\right|
 &\leq \sum_{r=1}^\infty \sum_{p\atop c\mid p}p_k(r)\left|(-1)^{ap+1}\sin\left(\frac{\pi a}{c}\right)\sum_h\frac{\omega_{h,p}^k}{\sin\left(\frac{\pi ah'}{c}\right)}e^{-\frac{\pi ia^2p_1h'}{c}-\frac{2\pi ihn}{p}+\frac{2\pi ih'}cm_{r,s}}\right|\\
 &\quad\times \left(\frac 2n\right)^{\frac k4-\frac 12}e^{2\pi+\frac{k\pi}{12}}e^{-\pi r}\int_{-\frac{1}{p(N+p)}}^{\frac{1}{p(N+p)}} d\Phi.
\end{align*}
Applying Lemma \ref{bound}, 
we may bound this by
\begin{align*}
  &\widehat{C}\sum_{r=1}^\infty p_k(r)e^{-\pi r}\sum_p p^{\varepsilon-\frac 12}\left(24n+k,p\right)^{\frac 12}n^{-\frac k4}
  \\
  &\leq \widehat{C}_1\sum_p p^{\varepsilon-\frac 12}\left(24n+k,p\right)^{\frac 12}n^{-\frac k4}
  \leq \widehat{C}_1\sum_{p\leq N} p^{\varepsilon-\frac 12}\sum_{d\mid p \atop d\mid 24n+k}d^{\frac 12}n^{-\frac k4}
  \\
  &\leq \widehat{C}_1 n^{-\frac k4}\sum_{d\mid 24n+k \atop d\leq N}d^{\frac 12}\sum_{p\leq N/d} (pd)^{\varepsilon-\frac 12}
  \leq \widehat{C}_2\sum_{d\mid 24n+k \atop d\leq N}d^{-\frac 12}N^{\varepsilon+\frac 12}n^{-\frac k4}
  \\
  &\leq \widehat{C}_2\sum_{d\mid 24n+k \atop d\leq N}d^{-\frac 12}n^{\frac{1-k}4+\frac{\varepsilon}2}=O(n^\varepsilon),
\end{align*}
where $\widehat{C},\widehat{C}_1$ and $\widehat{C}_2$ are constants. Here we bound trivially
\[\sum_{d|24n+k}d^{-\frac 12}<\sum_{d|24n+k}1=\sigma_0(24n+k)=o(n^\epsilon)\]
and choose $0<\epsilon<\varepsilon/2$, where $\sigma_0(n)$ denotes the number of the divisors of $n$ and satisfies $\sigma_0(n)=o(n^\epsilon)$ \cite[p.296]{Apostol-1976}.
Therefore $|S_{21}|=O(n^\varepsilon)$ for $k\leq12$.

Both $S_{22}$ and $S_{23}$ can be bounded in a similar manner, so we will focus on $S_{22}$. We can rewrite the integral as follows:
\[\int_{-\frac{1}{p(\widetilde{p}_1+p)}}^{-\frac{1}{p(N+p)}}=\sum_{\ell=\widetilde{p}_1+p}^{ N+p-1}\int_{-\frac{1}{p\ell}}^{-\frac{1}{p(\ell+1)}},\]
then
\begin{align*}
 \left|S_{22}\right|&\leq \left|\sum_{r=1}^\infty \sum_{p\atop c\mid p}p_k(r)\sum_{\ell=\widetilde{p}_1+p}^{ N+p-1}\int_{-\frac{1}{p\ell}}^{-\frac{1}{p(\ell+1)}} z^{\frac k2-1} e^{\frac {2\pi z}{p}\left(n-\frac {k}{24}\right)+\frac{k\pi}{12pz}}q_1^rd\Phi\right.
 \\
&\quad\left.\times(-1)^{ap+1}\sin\left(\frac{\pi a}{c}\right)\sum_h\frac{\omega_{h,p}^k}{\sin\left(\frac{\pi ah'}{c}\right)}e^{-\frac{\pi ia^2p_1h'}{c}-\frac{2\pi ihn}{p}+\frac{2\pi ih'}cm_{r,s}}\right|=:A_k.
\end{align*}
We use the condition $N<\widetilde{p}_1+p\leq\ell$ and so we can rearrange the summation from $\sum_{\ell=\widetilde{p}_1+p}^{ N+p-1}$
to $\sum_{\ell=N+1}^{N+p-1}$, but we also have to rewrite the sum over $h$ to count all the terms that contribute
\begin{align*}
 A_k&=\left|\sum_{r=1}^\infty \sum_{p\atop c\mid p}p_k(r)\sum_{\ell=N+1}^{N+p-1}\int_{-\frac{1}{p\ell}}^{-\frac{1}{p(\ell+1)}} z^{\frac k2-1} e^{\frac {2\pi z}{p}\left(n-\frac {k}{24}\right)+\frac{k\pi}{12pz}}q_1^rd\Phi\right.\\
&\quad\left.\times(-1)^{ap+1}\sin\left(\frac{\pi a}{c}\right)\sum_{h\atop N<\widetilde{p}_1+p\leq\ell}\frac{\omega_{h,p}^k}{\sin\left(\frac{\pi ah'}{c}\right)}e^{-\frac{\pi ia^2p_1h'}{c}-\frac{2\pi ihn}{p}+\frac{2\pi ih'}cm_{r,s}}\right|.
\end{align*}
Now by the theory of Farey fractions we have
\[\widetilde{p}_1\equiv -h'~~(\text{mod}~p),~\widetilde{p}_2\equiv h'~~(\text{mod}~p),~N-p\leq\widetilde{p}_i \leq N ,\]
for $i=1,2$. One can easily calculate that
\[\sum_{\ell=N+1}^{N+p-1}\int_{-\frac{1}{p\ell}}^{-\frac{1}{p(\ell+1)}}d\Phi\leq \frac{2}{p\sqrt{n}}.\]
Therefore, by applying Lemma \ref{bound}, we can bound $S_{22}$ in a manner similar to $S_{21}$. As a result, we obtain
\[S_{22}=O(n^\varepsilon),S_{23}=O(n^\varepsilon).\]
We now turn our attention to $S_1$.
The integral is split in the following way:
\[\int_{-\vartheta'_{h,p}}^{\vartheta''_{h,p}}
=\int_{-\frac{1}{pN}}^{\frac{1}{pN}}-\int_{-\frac{1}{pN}}^{-\frac{1}{p(\widetilde{p}_1+p)}}
-\int_{\frac{1}{p(\widetilde{p}_2+p)}}^{\frac{1}{pN}}.\]
We denote the corresponding sums by $S_{11}, S_{12},$ and $S_{13}$. It can be shown that $S_{12}$ and $S_{13}$ contribute to the error term. We begin with $S_{12}$. Following a similar analysis to the error terms of $S_2$, we write the integral as:
\[\int_{-\frac{1}{pN}}^{-\frac{1}{p(\widetilde{p}_1+p)}}=\sum_{\ell=N}^{ \widetilde{p}_1+p-1}\int_{-\frac{1}{p\ell}}^{-\frac{1}{p(\ell+1)}}.\]
Substituting into $S_{12}$ yields:
\begin{align*}
  S_{12}=\sum_{p\atop c\mid p}\sum_{\ell=N}^{ \widetilde{p}_1+p-1}\int_{-\frac{1}{p\ell}}^{-\frac{1}{p(\ell+1)}} z^{\frac k2-1}e^{\frac {2\pi z}{p}\left(n-\frac {k}{24}\right)+\frac{k\pi}{12pz}}d\Phi
\frac{i\sin\left(\frac{\pi a}{c}\right)}{(-1)^{ap+1}}\sum_{h}\frac{\omega_{h,p}^k}{\sin\left(\frac{\pi ah'}{c}\right)} e^{-\frac{\pi ia^2p_1h'}{c}-\frac{2\pi ihn}{p}}
\end{align*}
Now, given the condition $\widetilde{p}_1\leq N$, we have that $\ell\leq \widetilde{p}1+p-1\leq N+p-1$, which restricts the summation over $h$. We can bound $S_{12}$ by summing over more integrals:
\begin{align*}
  |S_{12}|&\leq\sum_{p\atop c\mid p}\sum_{\ell=N}^{N+p-1}\int_{-\frac{1}{p\ell}}^{-\frac{1}{p(\ell+1)}} \left|z^{\frac k2-1}e^{\frac {2\pi z}{p}\left(n-\frac {k}{24}\right)+\frac{k\pi}{12pz}}\right|d\Phi\\
&\quad\times\left|\sin\left(\frac{\pi a}{c}\right)(-1)^{ap+1}\sum_{h\atop \ell\leq \widetilde{p}_1+p-1\leq N+p-1}\frac{\omega_{h,p}^k}{\sin\left(\frac{\pi ah'}{c}\right)} e^{-\frac{\pi ia^2p_1h'}{c}-\frac{2\pi ihn}{p}}\right|\\
&=O(n^\varepsilon).
\end{align*}
Similarly, we have
\[S_{13}=O(n^\varepsilon).\]
Therefore,
\begin{align*}
  \Sigma_1=i\sin\left(\frac{\pi a}{c}\right)\sum_{h,p\atop c\mid p}\omega_{h,p}^k\frac{(-1)^{ap+1}}{\sin\left(\frac{\pi ah'}{c}\right)} e^{-\frac{\pi ia^2p_1h'}{c}-\frac{2\pi ihn}{p}}\int_{-\frac{1}{pN}}^{\frac{1}{pN}} z^{\frac k2-1}e^{\frac {2\pi z}{p}\left(n-\frac {k}{24}\right)+\frac{k\pi}{12pz}}d\Phi+O(n^\varepsilon).
\end{align*}
We now turn to $\Sigma_2$. Based on \eqref{C-k-4v}, we have
\begin{align}\label{C-k-v}
 &C_k\left(ah',\frac{lc}{c_1},c;q_1\right)
 \notag \\[5pt]
 &=\frac{i }{2(q_1)_{\infty}^k}\left(\sum_{m=0}^\infty\frac{(-1)^m e^{-\frac{\pi iah'}{c}}q_1^{\frac{m(m+1)}2+\frac l{2c_1}}}{1-e^{-\frac{2\pi iah'}{c}}q_1^{m+\frac l{c_1}}}-\sum_{m=1}^\infty\frac{(-1)^m e^{\frac{\pi iah'}{c}}q_1^{\frac{m(m+1)}2-\frac l{2c_1}}}{1-e^{\frac{2\pi iah'}{c}}q_1^{m-\frac l{c_1}}}\right)
 \notag \\[5pt]
 &=\frac{i }{2(q_1)_{\infty}^k}\left(\sum_{m=0}^\infty(-1)^m\sum_{r=0}^\infty e^{-\frac{\pi iah'}{c}-\frac{2\pi iah'r}{c}}q_1^{\frac{m(m+1)}2+\frac l{2c_1}+rm+\frac{rl}{c_1}}\right.
 \notag \\[5pt]
 &\left.\quad\quad\quad\quad\quad\quad-\sum_{m=1}^\infty(-1)^m\sum_{r=0}^\infty e^{\frac{\pi iah'}{c}+\frac{2\pi iah'r}{c}}q_1^{\frac{m(m+1)}2-\frac l{2c_1}+rm-\frac{rl}{c_1}}\right).
\end{align}
From this expression, we can deduce
\begin{align}\label{C-k-4v-t}
e^{-\frac{\pi ia^2p_1h'}{cc_1}+\frac{2\pi ialh'}{cc_1}+\frac{k\pi}{12pz}} q_1^{-\frac{l^2}{2c_1^2}} C_k\left(ah',\frac {lc}{c_1},c;q_1\right)
:=\sum_{r\geq r_0}\sum_{s~~(\text{mod}~c)}b_k(r,s)e^{\frac{2\pi i h'm_{a,c,p,r}}{p}}q_1^r.
\end{align}
We now clarify that $m_{a,c,p,r}\in\mathbb{Z}$ and $r_0$ may be negative. The part with negative $r$ contributes to the main part. We further rewrite \eqref{C-k-4v-t} using $1/(q_1)_{\infty}^k=1+O(q_1^k)$ inside of $C_k\left(ah',\frac {lc}{c_1},c;q_1\right)$. Consequently, the main contribution of
\[e^{-\frac{\pi ia^2p_1h'}{cc_1}+\frac{2\pi ialh'}{cc_1}+\frac{k\pi}{12pz}} q_1^{-\frac{l^2}{2c_1^2}} C_k\left(ah',\frac {lc}{c_1},c;q_1\right)\]
comes from
\begin{equation}\label{mian-contr}
  \pm\frac i2 e^{-\frac{\pi ia^2p_1h'}{cc_1}+\frac{2\pi ialh'}{cc_1}+\frac{k\pi}{12pz}} q_1^{-\frac{l^2}{2c_1^2}}(-1)^m q_1^{\frac{m(m+1)}2\pm\frac l{2c_1}+rm\pm\frac{rl}{c_1}}e^{\mp\frac{\pi iah'}{c}\mp\frac{2\pi iah'r}{c}}.
\end{equation}
From this, we can separate the expression into the roots of unity and the part that depends on the variable $z$. The roots of unity appear as follows
\[\exp\left(\frac{2\pi ih'}{p}\left(-\frac{a^2p_1p}{2cc_1}+\frac{lap}{cc_1}
-\frac{l^2}{2c_1^2}+rm\pm\frac{rl}{c_1}\mp\frac{pra}{c}
+\frac{m(m+1)}{2}\pm\frac{l}{2c_1}\mp\frac{ap}{2c}\right)\right).\]
Rewriting the expression in the second bracket and using the fact that $\frac{p_1}{c_1}=\frac pc$, along with rearranging the sum, it becomes possible to demonstrate that the contribution of the roots of unity looks like $\exp\left(\frac{2\pi ih'm_{a,c,p,r}}{p}\right)$, where $m_{a,c,p,r} \in\mathbb{Z}$. The intriguing part occurs for $\exp\left(\frac{\pi}{pz}T_k\right)$, where
\[T_k=\frac{l^2}{c_1^2}+\frac{k}{12}-2rm\mp2r\frac{l}{c_1}-m(m+1)\mp\frac{l}{c_1}.\]
This part contributes to the circle method exactly if $T_k>0$ which is equivalent to $-T_k<0$.
Firstly we treat the case with the plus sign in \eqref{mian-contr}. By multiplying by $(-1)$ and assuming $m>0$ we have
\[-T_k=-\frac{l^2}{c_1^2}-\frac{k}{12}+2rm+2r\frac{l}{c_1}+m(m+1)+\frac{l}{c_1}>0,\]
for $k\leq 12$. For $m=0$, define $r$ to be a solution to the following inequality:
\[ -\frac{l^2}{c_1^2}-\frac{k}{12}+2r\frac{l}{c_1}+\frac{l}{c_1}<0.\]
This is equivalent to $T_k>0$ and so this contributes to the main part in the Hardy-Ramanujan circle method.
Now choosing the minus sign in \eqref{mian-contr} that becomes
\[-T_k=-\frac{l^2}{c_1^2}-\frac{k}{12}+2rm-2r\frac{l}{c_1}+m(m+1)-\frac{l}{c_1},\]
which $>0$ for $m\geq2$ and $k\leq 12$. For $m=1$ we define $f:=[0,1]\rightarrow \mathbb{R}$ by
\[f(x)=:-x^2-x(1+2r)-\frac{k}{12}+2+2r.\]
Calculating the maximum and computing the values of the function we see that on the
boundary the function is negative, i.e., $f(1)=-\frac{k}{12}<0$. Thus this contributes to the
main part in the Circle Method. So there are two contributions coming from each of the
two terms of $C_k\left(ah',\frac {lc}{c_1},c;q_1\right)$. The first one comes from the first sum, if $m = 0$, and this
contributes with
\[\frac i2 e^{-\frac{\pi ia^2p_1h'}{cc_1}+\frac{2\pi ialh'}{cc_1}-\frac{\pi iah'}{c}+\frac{k\pi}{12pz}} q_1^{-\frac{l^2}{2c_1^2}+\frac l{2c_1}}\sum_{r\geq0\atop \delta_{a,c,p,r,k}^+>0 }e^{-\frac{2\pi iah'r}{c}}q_1^{\frac{rl}{c_1}},\]
where $\delta_{a,c,p,r,k}^+ = \frac{l^2}{2c_1^2}+\frac{k}{24}-(r+\frac 12)\frac{l}{c_1}$.The second contribution comes from the second
sum, if $m = 1$, and this contributes with
\[ \frac i2e^{-\frac{\pi ia^2p_1h'}{cc_1}+\frac{2\pi ialh'}{cc_1}+\frac{\pi iah'}{c}+\frac{k\pi}{12pz}} q_1^{-\frac{l^2}{2c_1^2}-\frac l{2c_1}+1}\sum_{r\geq0\atop \delta_{a,c,p,r,k}^->0 }e^{\frac{2\pi iah'r}{c}}q_1^{r\left(1-\frac{l}{c_1}\right)},\]
where $\delta_{a,c,p,r,k}^- = \frac{l^2}{2c_1^2}+\frac{k}{24}-1-r(1-\frac l{c_1})+\frac{l}{2c_1}$. So we get the main contributions of $\Sigma_2$:
\begin{align}\label{sigma2-main}
2\sin\left(\frac{\pi a}{c}\right)\sum_{\substack{p,r\\ c\nmid p\\ \delta_{a,c,p,r,k}^j>0\\ j\in\{+,-\} }}
(-1)^{ap+l}\sum_h\omega_{h,p}^k e^{\frac{2\pi i}{p}\left(-nh+m_{a,c,p,r}^jh'\right)}
\int_{-\vartheta'_{h,p}}^{\vartheta''_{h,p}} z^{\frac k2-1}e^{\frac {2\pi z}{p}\left(n-\frac {k}{24}\right)+\frac{2\pi}{pz}\delta_{a,c,p,r,k}^j} d\Phi.
\end{align}
Now it is possible to rewrite the sum over $p$ into the sum, where the $p$'s have the same
values for $c_1$ and $l$ and thus the $T_k$ are constant in each class and the condition
$T_k>0$ is independent of $p$ in each class. Moreover, it is clear as $c_1$ and $l$ are finite
numbers and for arbitrary large r there do not exist any solutions to $T_k>0$, so
that there are only finitely many solution to the inequality. That means it is possible to
split the sum over $r$ into positive $T_k$, which by the above argument is a finite sum
and into negative $T_k$, where the part with negative $T_k$ contributes to the error.
By symmetrizing the integral and now using Lemma \ref{bound}, it is possible to bound all
the terms exactly the same way we did for $\Sigma_1$:
\begin{align*}
 \Sigma_2=2\sin\left(\frac{\pi a}{c}\right)\sum_{\substack{p,r\\ c\nmid p\\ \delta_{a,c,p,r,k}^j>0 \\j\in\{+,-\} }}(-1)^{ap+l}
 &\sum_h\omega_{h,p}^k e^{\frac{2\pi i}{p}\left(-nh+m_{a,c,p,r}^jh'\right)}
 \\
&\times\int_{-\frac{1}{pN}}^{\frac{1}{pN}} z^{\frac k2-1}e^{\frac {2\pi z}{p}\left(n-\frac {k}{24}\right)+\frac{2\pi}{pz}\delta_{a,c,p,r,k}^j} d\Phi+O(n^\varepsilon).
\end{align*}
To finish the proof we have to evaluate integrals of the following form
\[I_{p,t}=\int_{-\frac{1}{pN}}^{\frac{1}{pN}}z^{\frac{k}{2}-1}e^{\frac{2\pi }{p}\left(z\left(n-\frac{k}{24}\right)+\frac{t}{z}\right)}d\Phi.\]
Upon substituting $z=\frac{p}{n}-ip\Phi$, we obtain
\begin{align}\label{hatS11-1}
I_{p,t}=\frac{1}{pi}\int_{\frac{p}{n}-\frac{i}{N}}^{\frac{p}{n}+\frac{i}{N}}z^{\frac{k}{2}-1}e^{\frac{2\pi }{p}\left(z\left(n-\frac{k}{24}\right)+\frac{t}{z}\right)}dz.
\end{align}
We now denote the circle through $\frac{p}{n}\pm\frac{i}{N}$ and tangent to the imaginary axis at $0$ by $\Gamma$. If $z=x+iy$, then $\Gamma$ is given by $x^2+y^2=\alpha x$, with $\alpha=\frac{p}{n}+\frac{n}{N^2p}$.
Using the fact that $2>\alpha>\frac{1}{p}$, $\mathrm{Re}(z)\leq\frac{p}{n}$ and $\mathrm{Re}\left(\frac{1}{z}\right)< p$ on the smaller arc, we can show that
the integral along the smaller arc is in $O\left(n^{-\frac{k}{8}}\right)$.
Moreover the path of integration in \eqref{hatS11-1} can be changed by Cauchy's Theorem into the larger arc of $\Gamma$. Thus
\[I_{p,t}=\frac{1}{pi}\int_{\Gamma}z^{\frac{k}{2}-1}e^{\frac{2\pi }{p}\left(z\left(n-\frac{k}{24}\right)+\frac{t}{z}\right)}dz+O(n^{-\frac{k}{8}}).\]
By transforming the circle to a straight line using $s=\frac{2\pi t}{pz}$, we derive
\[I_{p,t}=\frac{2\pi}{p}\left(\frac{2\pi t}{p}\right)^{\frac{k}{2}}\frac{1}{2\pi i}
\int_{\gamma-i\infty}^{\gamma+i\infty}
s^{-\frac{k}{2}-1}e^{s+\frac{\beta}{s}}ds+O(n^{-\frac{k}{8}}),\]
where $\gamma\in\mathbb{R}$ and $\beta=\frac{\pi^2 t}{6p^2}\left(24n-k\right)$.
According to the definition of the modified Bessel function of the first kind, we obtain
\[I_{p,t}=\frac{2\pi}{p}\left(\frac{24t}{24n-k}\right)^{\frac k4}I_{\frac k2}\left(\frac {\pi}p\sqrt{\frac{2t(24n-k)}{3}}\right)+O(n^{-\frac{k}{8}}).\]
Finally, we arrive at
\begin{align*}
  \Sigma_2+\Sigma_1
  =&2\sin\left(\frac{\pi a}{c}\right)\sum_{\substack{p,r\\ c\nmid p\\ \delta_{a,c,p,r,k}^j>0 \\j\in\{+,-\} }}D_{a,c,p,k}(-n,m_{a,c,p,r}^j)
\int_{-\frac{1}{pN}}^{\frac{1}{pN}} z^{\frac k2-1}e^{\frac {2\pi z}{p}\left(n-\frac {k}{24}\right)+\frac{2\pi}{pz}\delta_{a,c,p,r,k}^j} d\Phi
\\[5pt]
 &+i\sum_{p\atop c\mid p}B_{a,c,p,k}(-n,0)\int_{-\frac{1}{pN}}^{\frac{1}{pN}} z^{\frac k2-1}e^{\frac {2\pi z}{p}\left(n-\frac {k}{24}\right)+\frac{k\pi}{12pz}}d\Phi+O(n^\varepsilon).
\end{align*}
This concludes the proof of Proposition \ref{asym} after incorporating the expressions for $I_{p,t}$.
\qed

Utilizing Proposition \ref{asym}, we are able to derive the asymptotic formula of $M_k(a,c;n)$.

{\noindent\it{Proof of Theorem \ref{asym-M}.}}
Let $\zeta_c^\beta$ be defined as $e^{\frac{2\pi i\beta}{c}}$. By considering the orthogonality of the roots of unity, we can confirm that
\begin{align}\label{M}
 \sum_{n=0}^{\infty} M_k(a,c;n)q^n=\frac 1c\sum_{n=0}^{\infty}p_k(n)q^n+\frac 1c\sum_{\beta=1}^{c-1}\zeta_c^{-a\beta}C_k(\zeta_c^\beta;q).
\end{align}
Recall that Iskander, Jain, and Talvola \cite{Iskander-Jain-Talvola-2020} provided an exact formula for $p_k(n)$. Specifically, for $n > k/24$,
\begin{align}\label{p-k}
 p_k(n)=&2\pi \left(n-\frac{k}{24}\right)^{-\frac k4-\frac 12}\sum_{m=0}^{\lfloor\frac{k}{24}\rfloor}\left(\frac{k}{24}-m\right)^{\frac k4+\frac 12}p_k(m)
 \nonumber\\
 &\quad\times\sum_{p=1}^{\infty}\frac{A_{p,k}(n,m)}{p}I_{\frac k2+1}\left(\frac{4\pi}{p}\sqrt{\left(\frac{k}{24}-m\right)\left(n-\frac{k}{24}\right)}\right).
\end{align}
Here \[A_{p,k}(n,m):=\sum_{0\leq h<p \atop (h,p)=1}e^{\pi iks(h,p)+\frac{2\pi i}{p}(mH-nh)},\]
where $H$ denotes the inverse of $h$ modulo $p$ and $s(h, p)$ is the usual Dedekind sum.
Theorem \ref{asym-M} can easily be deduced from \eqref{M}, \eqref{p-k}, and Proposition \ref{asym}, thereby completing the proof.
\qed

Utilizing  Theorem \ref{asym-M}, we are able to demonstrate Theorem \ref{c-asym} through a concise analysis.

{\noindent\it{Proof of Theorem \ref{c-asym}.}}
In this case, we only need to show that the main contribution of \eqref{eq-asym-M} is the first term, which is indeed $\frac{p_k(n)}c$. Therefore, it is necessary to compare $p_k(n)$ with the Bessel function $I_{\frac k2}$.

Based on the works of \cite{Hardy-Ramanujan-1918-1, Rademacher-Zuckerman-1938}, we obtain
\begin{equation}\label{asym-pk}
  p_k(n)=2\left(\frac k3\right)^{\frac{1+k}{4}}(8n)^{-\frac{3+k}{4}}e^{\pi\sqrt{\frac{2kn}3}}\left(1+O\left(\frac{1}{\sqrt{n}}\right)\right).
\end{equation}
Moreover, the Bessel function is given by \cite[Eq. (4.12.7)]{Andrews-Askey-Roy-1999} as
\begin{equation*}
 I_\ell(x)=\frac{e^x}{\sqrt{2\pi x}}+O\left(\frac{e^x}{x^{\frac 32}}\right).
\end{equation*}
In the second term of \eqref{eq-asym-M}, given that $p\geq c$, we have
\[\frac{\pi\sqrt{k(24n-k)}}{6p}\leq\frac{\pi\sqrt{k(24n-k)}}{6c}=\pi\sqrt{\frac{2kn}{3c^2}-\frac{k^2}{36c^2}}<\pi\sqrt{\frac{2kn}{3}}. \]
Consequently, this term is significantly smaller than $\frac{p_k(n)}c$.

We assert that the same holds true for the contribution arising from the third term. To demonstrate this, observe that, directly from the definitions in \eqref{delta+} and \eqref{delta-}, we deduce $\delta_{a,c,p,r,k}^j\leq \frac{k}{24}$. Consequently,
\begin{equation*}
  \frac{\pi}{p}\sqrt{\frac{2 \delta_{\beta,c,p,r,k}^j(24n-k)}{3}}\leq \pi\sqrt{\frac{k(24n-k)}{36}}=\pi\sqrt{\frac{2kn}{3}-\frac{k^2}{36}}<\pi\sqrt{\frac{2kn}{3}}.
\end{equation*}
Combining \eqref{asym-pk}, we reach Theorem \ref{c-asym}.
\qed

\section{Proof of the inequalities for $M_k(a,c;n)$}\label{ineq}
In this section, we have provided bounds for the error terms discussed in the previous section, enabling the computation of Theorem \ref{ineq-M3}. Furthermore, we have directly demonstrated Theorems \ref{addi} and \ref{ineq-concave} by employing Theorem \ref{c-asym}.

Utilizing \eqref{M}, it can be readily deduced that
\begin{equation*}
  \sum_{n=0}^{\infty} \left(M_k(a,c;n)-M_k(b,c;n)\right)q^n=\frac 2c\sum_{\beta=1}^{\frac{c-1}2}\rho_\beta(a,b,c)C_k(\zeta_c^\beta;q),
\end{equation*}
where $\rho_\beta(a,b,c):=\cos\left(\frac{2\pi a\beta}{c}\right)-\cos\left(\frac{2\pi b\beta}{c}\right)$.
Based on Proposition \ref{asym}, we obtain
\begin{equation*}
  M_k(a,c;n)-M_k(b,c;n)=\sum_{\beta=1}^{\frac{c-1}2}\left(S_\beta(a,b,c;n)+\sum_{j\in\{+,-\}}T_\beta^j(a,b,c;n)+O(n^\varepsilon)\right),
\end{equation*}
where
\begin{align*}
S_\beta(a,b,c;n):=\rho_\beta(a,b,c)\frac{4\pi k^{\frac k4} i}{c\left(24n-k\right)^{\frac k4}}\sum_{1\leq p\leq N\atop c\mid p}\frac{B_{\beta,c,p,k}(-n,0)}{p}I_{\frac k2}\left(\frac{\pi\sqrt{k(24n-k)}}{6p}\right),
\end{align*}
and
\begin{align*}
T_\beta^j(a,b,c;n)&:=\rho_\beta(a,b,c)\frac{8\cdot 24^{\frac k4}\pi\sin\left(\frac{\pi \beta}{c}\right)}{c\left(24n-k\right)^{\frac k4}}
\sum_{\substack{1\leq p\leq N\\ c\nmid p \\ r\geq 0\\ \delta_{\beta,c,p,r,k}^j>0  }} \frac{\left(\delta_{\beta,c,p,r,k}^j\right)^{\frac k4}D_{\beta,c,p,k}(-n,m_{\beta,c,p,r}^j)}{p}
\\[5pt]
&\quad\quad\quad\quad\quad\quad\quad\quad\quad\quad\quad\quad \quad\quad\times I_{\frac k2}\left(\frac{\pi}{p}\sqrt{\frac{2 \delta_{\beta,c,p,r,k}^j(24n-k)}{3}}\right).
\end{align*}
From \cite[Eq.(1.6)]{Baricz-2008}, we have
\begin{equation}\label{def-I}
  I_{\frac 32}(x)=\sqrt{\frac{2}{\pi x}}\left(\cosh x-\frac{\sinh x}{x}\right),
\end{equation}
which is an increasing function for $x>0$. Additionally, we can confirm that the modified Bessel function of the first kind, $I_{\frac k2}(x)$, is increasing for $x>0$ and odd $k\leq12$ using Mathematica.

To identify the main contribution occurring in $ M_k(a,c;n)-M_k(b,c;n)$, we must determine the largest argument present in $S_\beta(a,b,c;n)$ and $T_\beta^j(a,b,c;n)$.
It is readily apparent that the largest argument of $S_\beta(a,b,c;n)$ and $T_\beta^j(a,b,c;n)$ is determined by the Bessel function $I_{\frac k2}$. Thus, we only need to identify the largest term arising from $I_{\frac k2}$. In $S_\beta$, the main term occurs for $p=c$, which is given by $I_{\frac k2}\left(\frac{\pi\sqrt{k(24n-k)}}{6c}\right)$.

In $T_\beta^j$, we first determine the largest term of $\delta_{\beta,c,1,r,k}^j$.
It is evident that the largest argument occurs when $r = 0$ for fixed $\beta$ and $p$.
Since $0<l<c$, we obtain $\delta_{\beta,c,p,0,k}^{-}<\delta_{\beta,c,p,0,k}^{+}$.
Consequently, we only need to examine the largest term of $\delta_{\beta,c,p,0,k}^+=\frac 12(\frac l c-\frac 12)^2-\frac 18+\frac{k}{24}$.
Assuming $\frac l c<\frac 12$, we recognize, by the symmetry of the parabola in the argument $\frac l c$, that $\delta_{\beta,c,p,0,k}^{+}$ attains its maximum for $l=1$, with the maximum value $\delta_0:=\frac{1}{2c^2}-\frac{1}{2c}+\frac{k}{24}$. Thus, if $l\neq 1$, we have $\delta_{\beta,c,p,r,k}^{j}< \delta_0$.
Furthermore,  it is clear that the largest term of $I_{\frac k2}$ in $T_\beta^j$ occurs if $r=0$, $l=1$, and $p=1$.
The last condition, combined with $l\equiv \beta p$ (mod $c$) and $0<l<c$, yields $l=\beta$. Hence, $D_{\beta,c,1,k}(-n,m_{\beta,c,1,0}^+)=1$.

Now, considering that
\begin{equation*}
  \sqrt{\frac{2\delta_0}3}>\frac{\sqrt{k}}{6c}
\end{equation*}
holds for $c\geq 2$ (respectively, $c\geq 4$) and odd integers $5\leq k\leq 12$ (respectively, $k=3$), the main contribution can be expressed as
\begin{align*}
  T_1:=\rho_1(a,b,c)\frac{8\cdot24^{\frac k4}\pi(\delta_0)^{\frac k4}\sin\left(\frac{\pi}{c}\right)}{c\left(24n-k\right)^{\frac k4}}I_{\frac k2}\left(\pi\sqrt{\frac{2 \delta_0(24n-k)}{3}}\right)
\end{align*}
and this value is positive for all $0\leq a<b\leq\frac {c-1}2$. Therefore, for sufficiently large $n$, we obtain $M_k(a,c;n)>M_k(b,c;n)$.
The following step entails providing a clear definition of `sufficiently large' by establishing bounds for the remaining terms.

\begin{lem}\label{b-main}
Let $\delta_0$ be defined as $\delta_0:=\frac{1}{2c^2}-\frac{1}{2c}+\frac{k}{24}$, as mentioned above. In this case, we obtain
\begin{align*}
|S_\beta(a,b,c;n)|\leq \frac{16k^{\frac k4}\left(1+\log(\frac{c-1}{2})\right)n^{\frac 12}}{c\left(24n-k\right)^{\frac k4}\left(1-\frac{\pi^2}{24}\right)}\cdot I_{\frac k2}\left(\frac{\pi\sqrt{k(24n-k)}}{6c}\right)=:\tilde{T}_0,
\end{align*}

\begin{align}\label{T-beta}
|T_\beta^j|\leq\frac{2k\cdot 24^{\frac k4}\pi(\delta_0)^{\frac k4}}{3\left(24n-k\right)^{\frac k4}} \left(n^{\frac 12}\cdot I_{\frac k2}\left(\pi\sqrt{\frac{\delta_0(24n-k)}{6}}\right)+I_{\frac k2}\left(\pi\sqrt{\frac{2\delta_0(24n-k)}{3}}\right)\right)=:\tilde{T}_1,
\end{align}
\begin{align*}
|T_1^+(a,b,c;n)-T_1|\leq\frac{2k\cdot 24^{\frac k4}\pi(\delta_0)^{\frac k4} n^{\frac 12}}{3\left(24n-k\right)^{\frac k4}} I_{\frac k2}\left(\pi\sqrt{\frac{\delta_0(24n-k)}{6}}\right)=:\tilde{T}_2,
\end{align*}
and
\begin{align*}
|T_1^-|\leq\frac{2k\cdot 24^{\frac k4}\pi(\delta_0)^{\frac k4} n^{\frac 12}}{3\left(24n-k\right)^{\frac k4}}I_{\frac k2}\left(\pi\sqrt{\frac{2\delta_0(24n-k)}{3}}\right)=:\tilde{T}_3,
\end{align*}
where $\beta\geq2$ in \eqref{T-beta}.

\end{lem}
\proof
Given that the largest argument in the Bessel function $I_{\frac{k}{2}}$ occurs when $p=c$ and $|\rho_j(a,b,c)|<2$, we can deduce the following
\begin{align}\label{S}
|S_\beta(a,b,c;n)|\leq \frac{8k^{\frac k4}\pi}{c\left(24n-k\right)^{\frac k4}}\cdot I_{\frac k2}\left(\frac{\pi\sqrt{k(24n-k)}}{6c}\right)\sum_{1\leq p\leq N\atop c\mid p}\frac{|B_{\beta,c,p,k}(-n,0)|}{p}.
\end{align}
Since $h$ and $h'$ traverse the same primitive residue classes modulo $p$, we modify the argument of the sine function in the summation from $\beta h' \rightarrow h$, and consequently shift to another representative of the equivalence class. With this adjustment, we can deduce
\begin{equation}\label{bound-B}
|B_{\beta,c,p,k}(-n,0)|\leq \sum_{h=1\atop (h,p)=1}^p
\frac 1{\left|\sin\left(\frac{\pi h}{c}\right)\right|}.
\end{equation}
Putting \eqref{bound-B} into \eqref{S}, we have
\begin{align*}
|S_\beta(a,b,c;n)|\leq \frac{8k^{\frac k4}\pi}{c\left(24n-k\right)^{\frac k4}}\cdot I_{\frac k2}\left(\frac{\pi\sqrt{k(24n-k)}}{6c}\right)\sum_{1\leq p\leq N\atop c\mid p}\frac{1}{p}\sum_{h=1\atop (h,p)=1}^p
\frac 1{\left|\sin\left(\frac{\pi h}{c}\right)\right|}.
\end{align*}
According to Ciolan's work \cite[\S4.2]{Ciolan-2022}, we can derive
\begin{equation}\label{sin}
 \sum_{h=1\atop (h,p)=1}^p\frac 1{\left|\sin\left(\frac{\pi h}{c}\right)\right|}\leq \frac{2p\left(1+\log(\frac{c-1}{2})\right)}{\pi\left(1-\frac{\pi^2}{24}\right)}.
\end{equation}
Consequently, we can deduce
\begin{align*}
|S_\beta(a,b,c;n)|\leq \frac{16k^{\frac k4}\left(1+\log(\frac{c-1}{2})\right)n^{\frac 12}}{c\left(24n-k\right)^{\frac k4}\left(1-\frac{\pi^2}{24}\right)}\cdot I_{\frac k2}\left(\frac{\pi\sqrt{k(24n-k)}}{6c}\right).
\end{align*}

Next, our aim is to establish a bound for $T_\beta^j$.
For $p\geq2$, we can trivially determine a bound as follows
\begin{equation*}
  |D_{\beta,c,p,k}(-n,m_{\beta,c,p,r}^j)|\leq p.
\end{equation*}
From the largest argument of $\delta_{\beta,c,p,r,k}^{j}$, we know that $\delta_{\beta,c,p,r,k}^{j}< \delta_0$. So we have
\begin{align*}
|T_\beta^j(a,b,c;n)|\leq\frac{16\cdot 24^{\frac k4}\pi(\delta_0)^{\frac k4}}{c\left(24n-k\right)^{\frac k4}} I_{\frac k2}\left(\pi\sqrt{\frac{\delta_0(24n-k)}{6}}\right)
\sum_{\substack{1\leq p\leq N\\ c\nmid p \\ r\geq 0\\ \delta_{\beta,c,p,r,k}^j>0 }}1
\end{align*}
for $\beta\geq2$ and $p\geq 2$.

We assert that the number of $r$ satisfying the condition $\delta_{\beta,c,p,r,k}^j>0$ can be bounded with respect to $c$. To begin with, let's examine the solution of $\delta_{\beta,c,p,r,k}^+>0$, which is given by $r<\frac{l}{2c}+\frac{kc}{24l}-\frac 12$. Let us define the function $f(l):=\frac{l}{2c}+\frac{kc}{24l}-\frac 12$, and consider it as a function of $l \in [1, c-1]$ for a fixed $c$. We calculate that the function attains its maximum for $l= 1$ or $l=c-1$, and this maximum can be bounded by $\frac{kc}{24}$. In other words, the number of $r$ satisfying the condition $\delta_{\beta,c,p,r,k}^+>0$ can be bounded by $\frac{kc}{24}$. Similarly, we obtain that the number of $r$ satisfying $\delta_{\beta,c,p,r,k}^->0$ can also be bounded by $\frac{kc}{24}$.
Hence, we can bound $T_\beta^j$ for $\beta\geq2$ and $p\geq 2$ using the following expression
\begin{align*}
\frac{2k\cdot 24^{\frac k4}\pi(\delta_0)^{\frac k4} n^{\frac 12}}{3\left(24n-k\right)^{\frac k4}} I_{\frac k2}\left(\pi\sqrt{\frac{\delta_0(24n-k)}{6}}\right).
\end{align*}
We bound the $p=1$ contribution by
\begin{align*}
\frac{2k\cdot 24^{\frac k4}\pi(\delta_0)^{\frac k4}}{3\left(24n-k\right)^{\frac k4}} I_{\frac k2}\left(\pi\sqrt{\frac{2\delta_0(24n-k)}{3}}\right).
\end{align*}
Summing the above two formulas together gives \eqref{T-beta}.
Repeating this step, we see that
\begin{align*}
|T_1^+(a,b,c;n)-T_1|&=\Bigg|\rho_1(a,b,c)\frac{8\cdot24^{\frac k4}\pi\sin\left(\frac{\pi}{c}\right)}{c\left(24n-k\right)^{\frac k4}}
\sum_{\substack{2\leq p\leq N\\ c\nmid p \\ r> 0\\ \delta_{1,c,p,r,k}^+>0  }} \frac{\left(\delta_{1,c,p,r,k}^j\right)^{\frac k4}D_{1,c,p,k}(-n,m_{1,c,p,r}^j)}{p}
\\[5pt]
&\quad\quad\times I_{\frac k2}\left(\frac{\pi\sqrt{2 \delta_{1,c,p,r,k}^j(24n-k)}}{3p}\right)\Bigg|
\\[5pt]
&\leq\frac{16\cdot 24^{\frac k4}\pi(\delta_0)^{\frac k4}}{c\left(24n-k\right)^{\frac k4}} I_{\frac k2}\left(\pi\sqrt{\frac{\delta_0(24n-k)}{6}}\right)
\sum_{\substack{2\leq p\leq N\\ c\nmid p \\ r> 0\\ \delta_{1,c,p,r,k}^+>0 }}1
\\[5pt]
&\leq\frac{2k\cdot 24^{\frac k4}\pi(\delta_0)^{\frac k4} n^{\frac 12}}{3\left(24n-k\right)^{\frac k4}} I_{\frac k2}\left(\pi\sqrt{\frac{\delta_0(24n-k)}{6}}\right).
\end{align*}
By the same analysis, we have
\begin{align*}
|T_1^-|\leq\frac{2k\cdot 24^{\frac k4}\pi(\delta_0)^{\frac k4} n^{\frac 12}}{3\left(24n-k\right)^{\frac k4}} I_{\frac k2}\left(\pi\sqrt{\frac{2\delta_0(24n-k)}{3}}\right),
\end{align*}
which completes the proof.
\qed

Next, we aim to provide an explicit bound for the $O(n^\varepsilon)$-term in Proposition \ref{asym}, taking into account the error terms originating from the circle method, the errors resulting from integrating over the remaining parts of the interval, and the errors introduced by integrating along the smaller arc.
From Section \ref{thm}, we are aware that $S_2$ contributes to the error of $\Sigma_1$ in the circle method. In $\Sigma_2$, this error term corresponds to the instances where $\delta_{\beta,c,p,r,k}^j$ is not positive, defined by $T_{err}$. Consequently, we present the following lemma to estimate the error terms arising from the circle method
\begin{lem}\label{b-error-1}
Define
\[c_1:=\sum_{n=1}^{\infty} p_k(n) e^{-\pi n},~~c_2:=\sum_{m=1}^{\infty} \frac{e^{-\pi m(m+1) / 2}}{1-e^{-\pi m}}~~\text{and}~~c_3:=\sum_{m=2}^\infty \frac{e^{-\frac{\pi m(m+1)}{2}}}{1-e^{-\pi m+\pi}}.\]
Let us denote $\delta_0$ by the expression $\frac{1}{2c^2} - \frac{1}{2c} + \frac{k}{24}$. For an odd integer $c$ and a positive integer $3\leq k \leq 12$, we obtain
\begin{align*}
   |S_2| \leq\frac{2^{\frac{k+6}4}e^{2\pi+\frac{k\pi}{12}}n^{\frac{2-k}4}\left(c_1+2\left(1+\left|\cos \left(\frac{\pi}{c}\right)\right|\right)(1+c_1)c_2\right)\left(1+\log(\frac{c-1}{2})\right)}{\pi\left(1-\frac{\pi^2}{24}\right)}=:\tilde{T}_4,
\end{align*}
and
\begin{align*}
|T_{err}|\leq 8\cdot 2^{\frac{k-2}4}e^{2\pi}n^{\frac{2-k}4}f(c)=:\tilde{T}_5,
\end{align*}
where
\begin{equation*}
  f(c):=\frac {e^{\frac{k\pi}{24}}(1+c_1e^{\pi\delta_0})}{\left(1-e^{-\frac{\pi }{c}}\right)}+\frac 12e^{\pi \delta_0+\frac{k\pi}{24}}(1+c_1)c_2+\frac 12 e^{\pi \delta_0+\frac{k\pi}{24}}(1+c_1)c_3.
\end{equation*}
\end{lem}
\proof
Using the fact that
\begin{align*}
 \operatorname{Re}(z)=\frac pn, ~~~ \frac p2\leq \operatorname{Re}\left(\frac 1z\right)<p,~~~ |z|^{\frac k2-1}\leq \left(\frac 2n\right)^{\frac k4-\frac 12}~ ~\text{and} ~~\vartheta_{h,p}'+\vartheta_{h,p}''\leq \frac 2{p\sqrt{n}},
\end{align*}
we have
\begin{align}\label{S-2-1}
   |S_2|&\leq 2^{\frac{k+2}4}e^{2\pi+\frac{k\pi}{12}}n^{-\frac k4}\sum_{1\leq p\leq N\atop c|p}\frac 1p\sum_{h=1\atop (h,p)=1}^{p}\frac 1 {\left|\sin\left(\frac{\pi h}{c}\right)\right|} \max_z \left|C_k\left(\zeta_c^h;q_1\right)-1\right|.
\end{align}
In this step, we apply \eqref{C-k-2}, \eqref{S2}, and a change of variables from $\beta h' \rightarrow h$. By incorporating \eqref{C-k}, we obtain
\begin{align*}
C_k\left(\zeta_c^h;q_1\right)
&=\frac{1}{\left(q_1\right)_{\infty}^k}+\frac{\left(1-\zeta_c^h\right)}{\left(q_1\right)_{\infty}^k} \sum_{m \in \mathbb{Z} \backslash\{0\}} \frac{(-1)^m q_1^{\frac{m(m+1)}{2}}}{1-\zeta_c^h q_1^m}
\\[5pt]
&=\frac{1}{\left(q_1\right)_{\infty}^k}+\frac{\left(1-\zeta_c^h\right)}{\left(q_1\right)_{\infty}^k} \sum_{m=1}^{\infty} \frac{(-1)^m q_1^{\frac{m(m+1)}{2}}}{1-\zeta_c^h q_1^m}+\frac{\left(1-\zeta_c^{-h}\right)}{\left(q_1\right)_{\infty}^k} \sum_{m=1}^{\infty} \frac{(-1)^m q_1^{\frac{m(m+1)}{2}}}{1-\zeta_c^{-h} q_1^m}.
\end{align*}
This follows that
\begin{align}\label{S-2-2}
\left|\left(C_k\left(\zeta_c^h;q_1\right)-1\right)\right|\leq & \sum_{n=1}^{\infty} p_k(n) e^{-\pi n}
\notag\\[5pt]
&+ \sum_{n=0}^{\infty} p_k(n) e^{-\pi n} \sum_{m=1}^{\infty} e^{-\pi m(m+1) / 2}\left|\frac{1-\zeta_c^h}{1-\zeta_c^h q_1^m}+\frac{1-\zeta_c^{-h}}{1-\zeta_c^{-h} q_1^m}\right| .
\end{align}
We bound the term further by noting that
\begin{align}\label{S-2-3}
\left|\frac{1-\zeta_c^h}{1-\zeta_c^h q_1^m}+\frac{1-\zeta_c^{-h}}{1-\zeta_c^{-h} q_1^m}\right| \leq 2 \frac{1+\left|\cos \left(\frac{\pi}{c}\right)\right|}{1-e^{-\pi m}} .
\end{align}
Combining \eqref{S-2-1}, \eqref{S-2-2} and \eqref{S-2-3}, we have
\begin{align*}
   |S_2|&\leq 2^{\frac{k+2}4}e^{2\pi+\frac{k\pi}{12}}n^{-\frac k4}\left(c_1+2\left(1+\left|\cos \left(\frac{\pi}{c}\right)\right|\right)(1+c_1)c_2\right)\sum_{1\leq p\leq N\atop c|p}\frac 1p\sum_{h=1\atop (h,p)=1}^{p}\frac 1 {\left|\sin\left(\frac{\pi h}{c}\right)\right|}
   \\[5pt]
   &\leq\frac{2^{\frac{k+6}4}e^{2\pi+\frac{k\pi}{12}}n^{\frac{2-k}4}\left(c_1+2\left(1+\left|\cos \left(\frac{\pi}{c}\right)\right|\right)(1+c_1)c_2\right)\left(1+\log(\frac{c-1}{2})\right)}{\pi\left(1-\frac{\pi^2}{24}\right)}.
\end{align*}
 In the last step we use \eqref{sin} again.
We continue by repeating this procedure for $T_{err}$.

We define $M_k(\beta h',l,c; q_1)$ to be the terms with positive exponents in the $q_1$-expansion of
\begin{equation*}
 e^{\frac{k\pi}{12pz}}q_1^{-\frac{l^2}{2c^2}}C_k\left(\beta h',l,c;q_1\right).
\end{equation*}
Doing the usual change of variable $\beta h'\rightarrow h$, we have
\begin{align*}
|T_{err}|\leq 8\cdot 2^{\frac{k-2}4}e^{2\pi}n^{-\frac k4}\sum_{h, p\atop c\nmid p}\frac 1p \max_z \left|M_k(h,l,c; q_1)\right|
\end{align*}
from the definition of $\Sigma_2$. By \eqref{C-k-v}, we break $C_k\left(h,l,c;q_1\right)$ down into four terms:
\begin{align*}
C_k\left(h,l,c;q_1\right)
=\frac{i}{2(q_1;q_1)_\infty^k}
\left(\sum_{r=0}^\infty e^{-\frac{\pi ih}{c}-\frac{2\pi irh}{c}}
q_1^{\frac{l}{2c}+\frac{rl}{c}}+\sum_{r=0}^\infty e^{\frac{\pi ih}{c}+\frac{2\pi irh}{c}}
q_1^{1-\frac{l}{2c}+r-\frac{rl}{c}}\right.
\\[5pt]
\left.\sum_{m=1}^\infty\frac{(-1)^m e^{-\frac{\pi ih}{c}}q_1^{\frac{m^2+m}{2}+\frac{l}{2c}}}{1-e^{-\frac{2\pi ih}{c}}q_1^{m+\frac lc}}-\sum_{m=2}^\infty\frac{(-1)^m e^{\frac{\pi ih}{c}}q_1^{\frac{m^2+m}{2}-\frac{l}{2c}}}{1-e^{\frac{2\pi ih}{c}}q_1^{m-\frac lc}}\right).
\end{align*}
Then we gain the following contribution to the first term of $M_k$:
\begin{align*}
   \left|\frac{e^{\frac{k\pi}{12}}}2 q_1^{-\frac{l^2}{2c^2}+\frac{l}{2c}}\sum_{m=0}^{\infty}p_k(m)q_1^m\sum_{r= r_0}^{\infty}q_1^{\frac{rl}{c}}\right|
  \leq \frac 12 e^{\frac{k\pi}{12}+\frac{\pi l^2}{2c^2}-\frac{\pi l}{2c}}\left( \sum_{r= r_0}^{\infty}e^{-\frac{\pi rl}{c}}+ \sum_{r=0}^{\infty}e^{-\frac{\pi rl}{c}}\sum_{m=1}^{\infty}p_k(m)e^{-\pi m}\right),
\end{align*}
where $r\geq r_0:=\lceil \frac{l}{2c}+\frac{kc}{24l}-\frac 12\rceil$ is a solution of $\delta_{\beta,c,p,r,k}^+<0$. Using
\begin{equation*}
  \sum_{r= r_0}^{\infty}e^{-\frac{\pi rl}{c}}=\frac{e^{-\frac{\pi r_0l}{c}}}{1-e^{-\frac{\pi l}{c}}},
\end{equation*}
and the usual geometric series we can bound the term further by
\begin{align*}
 \frac {e^{\frac{k\pi}{12}+\frac{\pi l^2}{2c^2}-\frac{\pi l}{2c}-\frac{\pi r_0l}{c}}}{2\left(1-e^{-\frac{\pi l}{c}}\right)}+\frac {e^{\frac{k\pi}{12}+\frac{\pi l^2}{2c^2}-\frac{\pi l}{2c}}c_1}{2\left(1-e^{-\frac{\pi l}{c}}\right)}
 &=\frac{e^{\frac{\pi l}{c}\left(\frac{kc}{12l}+\frac{l}{2c}-\frac 12-r_0\right)}\left(1+c_1e^{\frac{\pi r_0l}{c}}\right)}{2\left(1-e^{-\frac{\pi l}{c}}\right)}
 \\
 &\leq \frac {e^{\frac{k\pi}{24}}(1+c_1e^{\pi\delta_0})}{2\left(1-e^{-\frac{\pi }{c}}\right)}.
\end{align*}
We can bound the second term by the same way. In the third and fourth summand
all the terms will contribute to the error as shown in Proposition \ref{asym}. We have
\begin{align*}
  \left|\frac{i}{2}\cdot e^{\frac{k\pi}{12pz}}q_1^{-\frac{l^2}{2c^2}+\frac{l}{2c}}\frac{1}{(q_1;q_1)_\infty^k}\sum_{m=1}^\infty\frac{(-1)^m e^{-\frac{\pi ih}{c}}q_1^{\frac{m^2+m}{2}}}{1-e^{-\frac{2\pi ih}{c}}q_1^{m+\frac lc}}\right|
  &\leq \frac 12e^{\frac{k\pi}{12}+\frac{\pi l^2}{2c^2}-\frac{\pi l}{2c}}(1+c_1)\sum_{m=1}^\infty \frac{e^{-\frac{\pi m(m+1)}{2}}}{1-e^{-\pi m-\frac{\pi l}{c}}}\\
  &\leq \frac 12e^{\pi \delta_0+\frac{k\pi}{24}}(1+c_1)c_2.
\end{align*}
Proceeding by repeating this step for the fourth sum, we can get a bound for the fourth term of $M_k$:
\begin{equation*}
  \frac 12 e^{\pi \delta_0+\frac{k\pi}{24}}(1+c_1)c_3.
\end{equation*}
Based on the above estimates, we have
\begin{align*}
|T_{err}|\leq 8\cdot 2^{\frac{k-2}4}e^{2\pi}n^{\frac{2-k}4}f(c),
\end{align*}
which completes the proof.
\qed

Next, let us address the errors that emerge from integrating over the leftover portions of the interval. Previously, we employed
\[\int_{-\vartheta'_{h,p}}^{\vartheta''_{h,p}}
=\int_{-\frac{1}{pN}}^{\frac{1}{pN}}-\int_{-\frac{1}{pN}}^{-\frac{1}{p(\widetilde{p}_1+p)}}
-\int_{\frac{1}{p(\widetilde{p}_2+p)}}^{\frac{1}{pN}}.\]
By integrating \eqref{S1} and incorporating it into the first term of the main contribution, we obtain the following error term in $\Sigma_1$
\begin{align*}
S_{1err}
:=&-i\sin\left(\frac{\pi \beta}{c}\right)\sum_{h,p\atop c\mid p}\frac{\omega_{h,p}^k(-1)^{\beta p+1}}{\sin\left(\frac{\pi \beta h'}{c}\right)} e^{-\frac{\pi i\beta^2ph'}{c}-\frac{2\pi ihn}{p}}
\\
&\times\left(\int_{-\frac{1}{pN}}^{-\frac{1}{p(\widetilde{p}_1+p)}}
+\int_{\frac{1}{p(\widetilde{p}_2+p)}}^{\frac{1}{pN}}\right) z^{\frac k2-1}e^{\frac {2\pi z}{p}\left(n-\frac {k}{24}\right)+\frac{k\pi}{12pz}}d\Phi.
\end{align*}
Similarly, by \eqref{sigma2-main}, we have the error term in $\Sigma_2$
\begin{align*}
S_{2err}&:=-2\sin\left(\frac{\pi \beta}{c}\right)\sum_{\substack{p,r\\ c\nmid p\\ \delta_{\beta,c,p,r,k}^j>0 \\j\in\{+,-\} }}(-1)^{a\beta+l}\sum_h\omega_{h,p}^k e^{\frac{2\pi i}{p}\left(-nh+m_{\beta,c,p,r}^jh'\right)}
 \\
&\times\left(\int_{-\frac{1}{pN}}^{-\frac{1}{p(\widetilde{p}_1+p)}}
+\int_{\frac{1}{p(\widetilde{p}_2+p)}}^{\frac{1}{pN}}\right) z^{\frac k2-1}e^{\frac {2\pi z}{p}\left(n-\frac {k}{24}\right)+\frac{2\pi}{pz}\delta_{\beta,c,p,r,k}^j} d\Phi.
\end{align*}

Next, we present the bounds for $S_{1err}$ and $S_{2err}$ separately.
\begin{lem}\label{b-error-2}
Let $c$ be an odd integer and $k\leq 12$ be a positive integer. Then we have
\begin{align*}
|S_{1err}|\leq \frac{8\cdot 2^{\frac{k-2}4}e^{2\pi+\frac{k\pi}{12}}\left(1+\log\left(\frac{c-1}2\right)\right)n^{\frac{2-k}4}}{\pi(1-\pi^2/24)}=:\tilde{T}_6,
\end{align*}
and
\begin{align*}
|S_{2err}|\leq16\cdot 2^{\frac{k-2}4}e^{2\pi}n^{\frac{2-k}4}\frac {e^{2\pi \delta_0}}{1-e^{-\frac{2\pi}c}}=:\tilde{T}_7.
\end{align*}
\end{lem}
\proof
Analogous to the proof methods employed in Lemmas \ref{b-main} and \ref{b-error-1}, we can estimate
\begin{align*}
|S_{1err}|&\leq4\cdot 2^{\frac{k-2}4}e^{2\pi+\frac{k\pi}{12}}n^{-\frac k4}\sum_{1\leq p\leq \sqrt{n}\atop c|p }\frac 1p
\sum_{h=1\atop (h,p)=1}^{p}\frac 1{\left|\sin\left(\frac{\pi h}{c}\right)\right|}
\\
\quad&\leq \frac{8\cdot 2^{\frac{k-2}4}e^{2\pi+\frac{k\pi}{12}}\left(1+\log\left(\frac{c-1}2\right)\right)n^{-\frac k4}}{\pi(1-\pi^2/24)}\sum_{1\leq p\leq \sqrt{n}\atop c|p }1
\\
\quad&\leq \frac{8\cdot 2^{\frac{k-2}4}e^{2\pi+\frac{k\pi}{12}}\left(1+\log\left(\frac{c-1}2\right)\right)n^{\frac{2-k}4}}{\pi(1-\pi^2/24)}.
\end{align*}
In a manner entirely parallel to that of $S_{1err}$, it can be demonstrated
\begin{align*}
|S_{2err}|&\leq8\cdot 2^{\frac{k-2}4}e^{2\pi}n^{-\frac k4}\sum_{\begin{subarray}{c}
p, r\\
\delta_{\beta,c,p,r,k}^j>0\\
~j\in\{+,-\}
\end{subarray}}
e^{2\pi\delta_{\beta,c,p,r,k}^j}
\end{align*}
In the subsequent step, we assess the sum over $r$ for $j=+$ and provide a bound in terms of $\delta_0$. Since it represents the largest argument, we can similarly establish a bound for the term with $j=-$. Consequently, we focus on the case where $j=+$. The resulting sum over $r$ yields
\begin{align*}
|S_{2err}|
&\leq 16\cdot 2^{\frac{k-2}4}e^{2\pi}n^{-\frac k4}\sum_{p,r\atop\delta_{\beta,c,p,r,k}^+>0}e^{2\pi\delta_{\beta,c,p,r,k}^+}
\\
\quad &\leq 16\cdot 2^{\frac{k-2}4}e^{2\pi}n^{\frac{2-k}4}\sum_{r\leq r_0-1}e^{\frac {\pi l^2}{c^2}+\frac {k\pi}{12}-\frac {2\pi rl}c-\frac {\pi l}c}
\\
\quad &=16\cdot 2^{\frac{k-2}4}e^{2\pi}n^{\frac{2-k}4}\cdot\frac{e^{\frac {\pi l^2}{c^2}+\frac {k\pi}{12}-\frac {\pi l}c} \left(e^{-\frac{2\pi l r_0}c}-1\right)}{e^{-\frac{2\pi l}c}-1}
\\
\quad &\leq16\cdot 2^{\frac{k-2}4}e^{2\pi}n^{\frac{2-k}4}\frac {e^{2\pi \delta_0}}{1-e^{-\frac{2\pi}c}}.
\end{align*}
In this case, we utilized the fact that the summation over $r$ constitutes an error term when it begins with $r_0$. Thus, we must account for all the $r$-terms where $r \leq r_0 - 1$. This concludes the proof.
\qed

Lastly, we establish bounds for the errors resulting from integration along the smaller arc. We denote these errors as $\hat{S}_{1err}$, originating from $\Sigma_1$, and $\hat{S}_{2err}$, associated with $\Sigma_2$.
\begin{lem}\label{b-error-3}
Let $c$ be an odd integer and $k\leq 12$ be a positive integer, then
\begin{align*}
 |\hat{S}_{1err}|\leq\frac{4\left(\frac{2^{\frac{k+6}4}}{k+2}+\frac{2^{\frac{k+4}{4}}}k\right)\left(1+\log(\frac{c-1}{2})\right)e^{2 \pi+\frac{k\pi}{12} }n^{\frac{4-k}{8}}}{\pi\left(1-\frac{\pi^2}{24}\right)}=:\tilde{T}_8,
\end{align*}
and
\begin{align*}
  |\hat{S}_{2err}|\leq 8\left(\frac{2^{\frac{k+6}4}}{k+2}+\frac{2^{\frac{k+4}{4}}}k\right) \frac {e^{2 \pi+2\pi \delta_0}n^{\frac{4-k}{8}}}{1-e^{-\frac{2\pi}c}}=:\tilde{T}_9.
\end{align*}
\end{lem}
\proof
Recall that we have computed integrals of the following form in Section \ref{thm}:
$$
I_{p,t}=\frac{1}{pi}\int_{\frac{p}{n}-\frac{i}{N}}^{\frac{p}{n}+\frac{i}{N}}z^{\frac{k}{2}-1}e^{\frac{2\pi }{p}\left(z\left(n-\frac{k}{24}\right)+\frac{t}{z}\right)}dz.
$$
Now we denote the circle through $\frac{p}{n} \pm \frac{i}{N}$ and tangent to the imaginary axis at 0 by $\Gamma$. For $z=x+i y, \Gamma$ is given by $x^2+y^2=\frac{p}{n}+\frac{n}{N^2 p} x=: \alpha x$. The path of integration can be changed into the larger arc, while on the smaller arc we have the following bounds: $\frac{1}{p}<\alpha<2, \operatorname{Re}(z) =\frac{p}{n}$ and $\operatorname{Re}\left(z^{-1}\right)<p$. This can be used to bound the integral over the smaller arc which we denoted by $\Gamma_S$. Splitting $I_{p, t}=I_{p, t}^{main }+I_{p, t}^{err }$ we can bound $I_{p, t}^ {err }$
$$
\begin{aligned}
|I_{p, t}^{e r r}| & \leq \frac{2}{p} e^{2 \pi+2 \pi t} \int_{\Gamma_S}|z|^{\frac{k}{2}-1} d z
\\[5pt]
&\leq \frac{2}{p} e^{2 \pi+2 \pi t}\left|\int_0^{\frac{p}{n}}\left(x^2+y^2\right)^{\frac{k-2}{4}}(d x+i d y)\right|
\\[5pt]
&=\frac{2}{p} e^{2 \pi+2 \pi t} \alpha^{\frac{k-2}{4}}\left|\int_0^{\frac{p}{n}} x^{\frac{k-2}{4}} d x+i \int_0^{\frac{p}{n}} x^{\frac{k-2}{4}} d y\right|
\\[5pt]
&=\frac{2}{p} e^{2 \pi+2 \pi t} \alpha^{\frac{k-2}{4}}\left|\frac{4}{k+2}\left(\frac{p}{n}\right)^{\frac{k+2}{4}}+i \int_0^{\frac{p}{n}} x^{\frac{k-2}{4}} \frac{d y}{d x} d x\right|
\\[5pt]
&=\frac{2}{p} e^{2 \pi+2 \pi t} \alpha^{\frac{k-2}{4}}\left|\frac{4}{k+2}\left(\frac{p}{n}\right)^{\frac{k+2}{4}}+i \int_0^{\frac{p}{n}} x^{\frac{k-4}{4}} \frac{\alpha-2 x}{2 \sqrt{\alpha-x}} d x\right| .
\end{aligned}
$$
 Define $h(x):=\frac{\alpha-2 x}{2 \sqrt{\alpha-x}}$. By computing the derivative of $h$, we see that the function $h(x)$ has its maximum at $x=0$. So we can bound further:
$$
\begin{aligned}
|I_{p, t}^{e r r} |& \leq \frac{2}{p} e^{2 \pi+2 \pi t} \alpha^{\frac{k-2}{4}}\left|\frac{4}{k+2}\left(\frac{p}{n}\right)^{\frac{k+2}{4}}+ \frac{i\alpha^{\frac{1}{2}}} {2} \int_0^{\frac{p}{n}} x^{\frac{k-4}{4}} d x\right|
\\[5pt]
& \leq \frac{2}{p} e^{2 \pi+2 \pi t}\left(\frac{4}{k+2}\left(\frac{p}{n}\right)^{\frac{k+2}{4}} \alpha^{\frac{k-2}{4}}+\frac{2}{k}\alpha^{\frac{k}{4}}\left(\frac{p}{n}\right)^{\frac{k}{4}}\right)
\\[5pt]
& \leq \frac{2}{p} e^{2 \pi+2 \pi t}\left(\frac{2^{\frac{k+6}4}}{k+2}+\frac{2^{\frac{k+4}{4}}}k\right) n^{-\frac{k}{8}}.
\end{aligned}
$$
Here it is used that $\frac{p}{n} \leq n^{-1 / 2}$ and $\alpha<2$. Combining the contributions from $\Sigma_1$, using usual formulas like \eqref{sin}, and estimation of the sum over $p$, we can bound the whole contribution by
\begin{align*}
 |\hat{S}_{1err}|\leq\frac{4\left(\frac{2^{\frac{k+6}4}}{k+2}+\frac{2^{\frac{k+4}{4}}}k\right)\left(1+\log(\frac{c-1}{2})\right)e^{2 \pi+\frac{k\pi}{12} }n^{\frac{4-k}{8}}}{\pi\left(1-\frac{\pi^2}{24}\right)}.
\end{align*}
The same can be done for $\Sigma_2$, we have
\begin{align*}
  |\hat{S}_{2err}|\leq&8\left(\frac{2^{\frac{k+6}4}}{k+2}+\frac{2^{\frac{k+4}{4}}}k\right) e^{2 \pi} n^{-\frac{k}{8}} \sum_{p,r\atop\delta_{\beta,c,p,r,k}^+>0}e^{2\pi\delta_{\beta,c,p,r,k}^+}
  \\
  &\leq 8\left(\frac{2^{\frac{k+6}4}}{k+2}+\frac{2^{\frac{k+4}{4}}}k\right) \frac {e^{2 \pi+2\pi \delta_0}n^{\frac{4-k}{8}}}{1-e^{-\frac{2\pi}c}}.
\end{align*}
This finishes the proof.
\qed

We will now proceed to prove Theorems \ref{ineq-M3}-\ref{ineq-concave}.

{\noindent\it{Proof of Theorem \ref{ineq-M3}.}}
By integrating the estimation of the main term and incorporating Lemmas \ref{b-main}-\ref{b-error-3}, we can deduce that
\begin{equation*}
  N_{a,b,c,k}=\min\left\{n\in\mathbb{N}\Bigg|T_1-\sum_j\tilde{T}_j>0\right\}.
\end{equation*}
This brings the proof of Theorem \ref{ineq-M3} to a conclusion.
\qed

{\noindent\it{Proof of Theorem \ref{addi}.}}
Applying Theorem \ref{c-asym}, we obtain
\begin{align*}
 \frac{M_k(a,c;n_1)M_k(a,c;n_2)}{M_k(a,c;n_1+n_2)}
 &\sim \frac{2}{8^{\frac{3+k}{4}}c}\left(\frac k3\right)^{\frac{1+k}{4}}\cdot \left(\frac{n_1+n_2}{n_1n_2}\right)^{\frac{3+k}{4}}\cdot\frac{e^{\pi\sqrt{\frac{2k}3}\left(\sqrt{n_1}+\sqrt{n_2}\right)}}
 {e^{\pi\sqrt{\frac{2k}3}\sqrt{n_1+n_2}}},
\end{align*}
which exceeds $1$ when $n_1$ and $n_2$ are sufficiently large.
This concludes the proof.
\qed

{\noindent\it{Proof of Theorem \ref{ineq-concave}.}}
Analogous to the proof of Theorem \ref{addi}, we examine the ratio
\begin{equation*}
 \frac{M_k(a,c;n_1)M_k(a,c;n_2)}{M_k(a,c;n_1-1)M_k(a,c;n_2+1)}.
\end{equation*}
Taking into account Theorem \ref{c-asym}, this ratio approaches
\begin{align*}
 \left(\frac{(n_1-1)(n_2+1)}{n_1n_2}\right)^{\frac{3+k}{4}}\cdot\frac{e^{\pi\sqrt{\frac{2k}3}\left(\sqrt{n_1}+\sqrt{n_2}\right)}}
 {e^{\pi\sqrt{\frac{2k}3}\left(\sqrt{n_1-1}+\sqrt{n_2+1}\right)}}>1
\end{align*}
when $n_1$ and $n_2$ are sufficiently large, with $n_1<n_2+1$.
This concludes the proof.
\qed

\section{Strict log-subadditivity for $k=3$}\label{SLS-3}
In this section, our objective is to prove the log-subadditivity property for $M_3(a,c;n)$ by providing explicit bounds for the error terms.
Let $\mu(n):=\pi\sqrt{2n-\frac 14}$. Then, according to Theorem \ref{asym-M}, we obtain
\begin{align}\label{M3-formula}
M_3(a,c;n)=
&\frac{p_3(n)}{c}+\frac{1}{c}\sum_{\beta=1}^{c-1}\zeta_c^{-a\beta}\frac{2\pi i}{\left(8n-1\right)^{\frac 34}}\sum_{1\leq p\leq N\atop c\mid p}\frac{B_{\beta,c,p,3}(-n,0)}{p}I_{\frac 32}\left(\frac{\mu(n)}{p}\right)
\notag\\[5pt]
&+\frac{1}{c}\sum_{\beta=1}^{c-1}\zeta_c^{-a\beta}\frac{4\pi\sin\left(\frac{\pi \beta}{c}\right)}{\left(n-\frac 18\right)^{\frac 34}}
\sum_{\substack{1\leq p\leq N\\ c\nmid p \\ r\geq 0\\ \delta_{\beta,c,p,r,3}^j>0  }} \frac{\left(\delta_{\beta,c,p,r,3}^j\right)^{\frac 34}D_{\beta,c,p,3}(-n,m_{\beta,c,p,r}^j)}{p}
\notag\\[5pt]
&\times I_{\frac 32}\left(\frac{2\sqrt{2 \delta_{\beta,c,p,r,3}^j}\mu(n)}{p}\right)+O(n^\varepsilon)
\notag\\[5pt]
=:&\frac{p_3(n)}{c}+G_1+G_2+O(n^\varepsilon).
\end{align}
To complete the proof of Theorem \ref{ineq-sum}, we require both an upper bound and a lower bound for $M_3(a,c;n)$. We begin by establishing bounds for $G_1$ and $G_2$.
\begin{lem}\label{lem-bound-G}
Letting $c$ be an odd integer and $n$ be a positive integer, we have
\begin{equation*}
  |G_1|\leq\frac{1.0816c^{\frac 32}\cdot n^{\frac 12}}{8n-1}\cdot e^{\frac{\mu(n)}{c}},
\end{equation*}
and
\begin{equation*}
  |G_2|\leq \frac{\frac c2\cdot n^{\frac 12}}{8n-1}\cdot e^{2\sqrt{2\delta_0}\mu(n)}.
\end{equation*}
\end{lem}
\pf
By \eqref{def-I}, we have
\begin{align}\label{b-I}
 I_{\frac 32}(x)
 &=\sqrt{\frac{2}{\pi x}}\left(\frac{e^x+e^{-x}}{2}-\frac{e^x-e^{-x}}{2x}\right)
 \notag\\[5pt]
 &=\frac{1}{\sqrt{2\pi x}}\left(\left(1-\frac 1x\right)e^x+\left(1+\frac 1x\right)e^{-x}\right)
 \notag\\[5pt]
 &\leq\frac{e^x}{\sqrt{2\pi x}}.
\end{align}
Incorporating both \eqref{bound-B} and \eqref{sin}, we deduce
\begin{equation*}
  |B_{\beta,c,p,3}(-n,0)|\leq \frac{2p\left(1+\log(\frac{c-1}{2})\right)}{\pi\left(1-\frac{\pi^2}{24}\right)}.
\end{equation*}
Combining this result with \eqref{b-I}, we obtain
\begin{align*}
  |G_1|
  &\leq \frac{2\pi }{\left(8n-1\right)^{\frac 34}}\sum_{1\leq p\leq N\atop c\mid p} \frac{2\left(1+\log(\frac{c-1}{2})\right)}{\pi\left(1-\frac{\pi^2}{24}\right)}I_{\frac 32}\left(\frac{\mu(n)}{p}\right)
  \\[5pt]
  &\leq  \frac{4\left(1+\log(\frac{c-1}{2})\right)\cdot n^{\frac 12}}{\left(8n-1\right)^{\frac 34}\left(1-\frac{\pi^2}{24}\right)}I_{\frac 32}\left(\frac{\mu(n)}{c}\right)
  \\[5pt]
  &\leq \frac{4\left(1+\log(\frac{c-1}{2})\right)\cdot n^{\frac 12}\sqrt{c}}{\pi(8n-1)\left(1-\frac{\pi^2}{24}\right)}\cdot e^{\frac{\mu(n)}{c}}
   \\[5pt]
  &\leq \frac{1.0816c^{\frac 32}\cdot n^{\frac 12}}{8n-1}\cdot e^{\frac{\mu(n)}{c}},
\end{align*}
where in the last step, we employ \cite[Eq.(3.13)]{Zhang-Zhong-2023}, which states that
\begin{equation}\label{b-log}
  \frac{1+\log(\frac{c-1}{2})}{\pi\left(1-\frac{\pi^2}{24}\right)}<0.2704c.
\end{equation}
Moreover, by the argument of main term in Section \ref{ineq}, we know that $\delta_{\beta,c,p,r,3}^{j}< \delta_0:=\frac{1}{2c^2}-\frac{1}{2c}+\frac 18\leq\frac 18$. Combining $|D_{\beta,c,p,3}(-n,m_{\beta,c,p,r}^j)\leq p$, we have
\begin{align*}
  |G_2|\leq \frac{4\pi}{\left(8n-1\right)^{\frac 34}} I_{\frac 32}\left(2\sqrt{2 \delta_0}\mu(n)\right) \sum_{\substack{1\leq p\leq N\\ c\nmid p \\ r\geq 0\\ \delta_{\beta,c,p,r,3}^j>0  }}.
\end{align*}
From the proof of Lemma \ref{b-main}, we obtain
\begin{equation*}
\sum_{\substack{ r\geq 0\\ \delta_{\beta,c,p,r,3}^j>0  }}\leq \frac c8.
\end{equation*}
Consequently, we deduce that
\begin{align*}
  |G_2|
  &\leq \frac{\pi cn^{\frac 12}}{2\left(8n-1\right)^{\frac 34}} I_{\frac 32}\left(2\sqrt{2 \delta_0}\mu(n)\right)
  \\
  &\leq \frac{\frac c2\cdot n^{\frac 12}}{8n-1}\cdot e^{2\sqrt{2\delta_0}\mu(n)}.
\end{align*}
Here, we use \eqref{b-I} once more, which concludes the proof.
\qed

We next give a bound for the error term $O(n^\varepsilon)$.
\begin{lem}\label{lem-bound-error}
Letting $c$ be an odd integer and $n$ be a positive integer, the error term $O(n^\varepsilon)$ can be bounded by
\[(392785.264c+178.307)n^{\frac 18}.\]
\end{lem}
\pf
By the definitions of constants $c_1,c_2$ and $c_3$ in Lemma \ref{b-error-1}, we can compute that $c_2<0.0085$ and $c_3<0.00002$.
Utilizing \eqref{asym-pk}, we deduce that for $n\geq1$
\begin{equation*}
  p_3(n)\sim \frac{2}{(8n)^{\frac 32}}e^{\pi\sqrt{2n}}\ll e^{\pi\sqrt{2n}},
\end{equation*}
which yields $c_1<2.7439$.
Therefore, by combining $\frac{1}{1-e^{-\frac{\pi}{c}}}\leq \pi c$ and $\delta_0\leq\frac 18$, the function $f(c)$ in Lemma \ref{b-error-1} can be bounded by
\[23.56c+0.035.\]
Consequently, when combined with Lemma \ref{b-error-1} and \eqref{b-log}, we obtain
\begin{align*}
 |S_2|\leq4337.44cn^{-\frac 14},
\end{align*}
and
\begin{align*}
 |T_{err}|\leq(120026.044c+178.307)n^{-\frac 14}.
\end{align*}
Considering Lemma \ref{b-error-2} and \eqref{b-log}, we deduce that
\begin{align*}
 |S_{1err}|\leq3021.35cn^{-\frac 14},
\end{align*}
and
\begin{align*}
 |S_{2err}|\leq140411.98cn^{-\frac 14}.
\end{align*}
Similarly, using Lemma \ref{b-error-3}, we obtain
\begin{align*}
 |\hat{S}_{1err}|\leq 2632.82cn^{\frac 18},
\end{align*}
and
\begin{align*}
 |\hat{S}_{2err}|\leq 122355.63cn^{\frac 18}.
\end{align*}
Summing up the six bounds above yields
\[O(n^\varepsilon)\leq(267796.814c+178.307)n^{-\frac 14}+124988.45cn^{\frac 18}.\]
Finally, taking the sum of the two coefficients and the highest power of $n$ from the two terms completes the proof.
\qed

Next, we present the upper and lower bounds for $p_3(n)$, which play a crucial role in proving Theorem \ref{ineq-sum}.
\begin{lem}\label{lem-bound-p3}
For $n\geq 1$, we have
\begin{equation}\label{bound-p3}
  \frac{2}{(8n-1)^{\frac 98}}\left(1-\frac{3}{\mu(n)}\right)e^{\mu(n)}<p_3(n)<\frac{2}{(8n-1)^{\frac 98}}\left(\frac 97+\frac{3}{\mu(n)}\right)e^{\mu(n)}.
\end{equation}
\end{lem}
\pf
Observe that, according to \cite[Eq.~9.6.26]{Abramowitz-Stegun-1964},
\begin{equation}\label{ditui}
  I_{v-1}(x)-I_{v+1}(x)=\frac{2v}{x}I_{v}(x).
\end{equation}
Furthermore, the modified Bessel function of the first kind with order $\frac 12$ is given by \cite[Eq.(1.5)]{Baricz-2008} as
\begin{equation}\label{def-I-12}
  I_{\frac 12}(x)=\sqrt{\frac{2}{\pi x}}\sinh x.
\end{equation}
By setting $v=\frac 32$ in \eqref{ditui} and combining it with \eqref{def-I} and \eqref{def-I-12}, we obtain
\begin{equation*}
  I_{\frac 52}(x)=\frac{1}{\sqrt{2\pi x}}\left(-\frac{6}{x}\cosh(x)+\left(2+\frac{6}{x^2}\right)\sinh(x)\right).
\end{equation*}
Hence, using \eqref{p-k}, we deduce that
\begin{align}\label{p3-f}
p_3(n)
&=\frac{2\pi}{(8n-1)^{\frac 78}}\sum_{k=1}^\infty \frac{A_{k,3(n,0)}}{k}I_{\frac 52}\left(\frac{\mu(n)}{k}\right)
\notag\\
&=\frac{2}{(8n-1)^{\frac 98}}\sum_{k=1}^N\frac{A_{k,3}(n,0)}{\sqrt{k}}\left[\left(\frac{3k^2}{\mu(n)^2}-\frac{3k}{\mu(n)}+1\right)e^{\frac{\mu(n)}{k}}\right.
\notag\\
&\quad\quad\quad\quad\quad\quad\quad\quad\quad\quad\quad\quad\left.-\left(\frac{3k^2}{\mu(n)^2}+\frac{3k}{\mu(n)}+1\right)e^{-\frac{\mu(n)}{k}}\right]+R(n,N),
\end{align}
where
\begin{align*}
R(n,N):=\frac{2}{(8n-1)^{\frac 98}}\sum_{k=N+1}^\infty\frac{A_{k,3}(n,0)}{\sqrt{k}}
&\left(-\frac{6k}{\mu(n)}\cosh\left(\frac{\mu(n)}{k}\right)\right.
\\
&\left.+\left(2+\frac{6k^2}{\mu(n)^2}\right)\sinh\left(\frac{\mu(n)}{k}\right)\right).
\end{align*}
Letting
\[f_n(k):=-\frac{6k}{\mu(n)}\cosh\left(\frac{\mu(n)}{k}\right)+
\left(2+\frac{6k^2}{\mu(n)^2}\right)\sinh\left(\frac{\mu(n)}{k}\right),\]
by Taylor's Theorem, we have
\begin{align*}
f_n(k)
&=\left(2+\frac{6k^2}{\mu(n)^2}\right)\sum_{m=0}^\infty\frac{\left(\frac{\mu(n)}{k}\right)^{2m+1}}{(2m+1)!}
-\frac{6k}{\mu(n)}\sum_{m=0}^\infty\frac{\left(\frac{\mu(n)}{k}\right)^{2m}}{(2m)!}
\\[5pt]
&=2\sum_{m=0}^\infty\frac{\left(\frac{\mu(n)}{k}\right)^{2m+1}}{(2m+1)!}+6\sum_{m=0}^\infty\frac{\left(\frac{\mu(n)}{k}\right)^{2m-1}}{(2m+1)!}
-6\sum_{m=0}^\infty\frac{\left(\frac{\mu(n)}{k}\right)^{2m-1}}{(2m)!}
\\[5pt]
&=2\sum_{m=0}^\infty\frac{\left(\frac{\mu(n)}{k}\right)^{2m+1}}{(2m+1)!}-6\sum_{m=1}^\infty\frac{2m}{(2m+1)!}\left(\frac{\mu(n)}{k}\right)^{2m-1}
\\[5pt]
&=\sum_{m=1}^\infty\frac{8m(m+1)}{(2m+3)!}\left(\frac{\mu(n)}{k}\right)^{2m+1}.
\end{align*}
It is obvious that $|A_{k,3}(n,0)|\leq k$. Since $\sum_{k=N+1}^\infty\sqrt{k}\cdot\frac{\left(\frac{\mu(n)}{k}\right)^{2m+1}}{(2m+1)!}$ is convergent, we have
\begin{align*}
 \left| \sum_{k=N+1}^\infty\frac{A_{k,3}(n,0)}{\sqrt{k}}f_n(k)\right|
 &\leq \sum_{k=N+1}^\infty\sqrt{k}f_n(k)
 \\[5pt]
 &=\sum_{k=N+1}^\infty\sum_{m=1}^\infty\frac{8m(m+1)\sqrt{k}}{(2m+3)!}\left(\frac{\mu(n)}{k}\right)^{2m+1}
 \\[5pt]
 &=\sum_{m=1}^\infty\frac{8m(m+1)\mu(n)^{2m+1}}{(2m+3)!}\sum_{k=N+1}^\infty k^{-2m-\frac 12}
 \\[5pt]
 &\leq\sum_{m=1}^\infty\frac{8m(m+1)\mu(n)^{2m+1}}{(2m+3)!}\int_N^\infty x^{-2m-\frac 12}dx
 \\[5pt]
 &=\sum_{m=1}^\infty\frac{8mN^{\frac 32}}{(2m+3)(4m-1)}\cdot\frac{\left(\frac{\mu(n)}{N}\right)^{2m+1}}{(2m+1)!}.
\end{align*}
Hence
\begin{align}\label{b-R(n,N)}
  |R(n,N)|
  &\leq\frac{8N^{\frac 32}}{7(8n-1)^{\frac 98}}\sum_{m=1}^\infty\frac{14m}{(2m+3)(4m-1)}\cdot\frac{\left(\frac{\mu(n)}{N}\right)^{2m+1}}{(2m+1)!}
  \notag\\
  &\leq\frac{8N^{\frac 32}}{7(8n-1)^{\frac 98}}\sum_{m=0}^\infty\frac{\left(\frac{\mu(n)}{N}\right)^{2m+1}}{(2m+1)!}
  \notag\\
  &=\frac{8N^{\frac 32}}{7(8n-1)^{\frac 98}}\sinh\left(\frac{\mu(n)}{N}\right).
\end{align}
Here we use that $(2m+3)(4m-1)>14m$ for $m\geq1$.

Let $N=1$ in \eqref{p3-f}, we have
\begin{align*}
p_3(n)=\frac{2}{(8n-1)^{\frac 98}}\left[\left(\frac{3}{\mu(n)^2}-\frac{3}{\mu(n)}+1\right)e^{\mu(n)}
-\left(\frac{3}{\mu(n)^2}+\frac{3}{\mu(n)}+1\right)e^{-\mu(n)}\right]+R(n,1).
\end{align*}
This together with \eqref{b-R(n,N)} gives
\begin{align*}
p_3(n)
&\leq\frac{2}{(8n-1)^{\frac 98}}\left[\left(\frac{3}{\mu(n)^2}-\frac{3}{\mu(n)}+\frac 97\right)e^{\mu(n)}
-\left(\frac{3}{\mu(n)^2}+\frac{3}{\mu(n)}+\frac 97\right)e^{-\mu(n)}\right]
\\
&<\frac{2}{(8n-1)^{\frac 98}}\left(\frac{3}{\mu(n)}+\frac 97\right)e^{\mu(n)}.
\end{align*}
The left-hand side of \eqref{bound-p3} can be easily verified. With this, we complete the proof.
\qed

Based on the above estimations, we now can bound $M_3(a,c;n)$ by $p_3(n)$.
\begin{prop}
Let $c$ be an odd integer. For $c\geq4$, we have
\begin{equation}\label{b-M-4}
  \frac{1}{2c}p_3(n)<M_3(a,c;n)< \frac{3}{2c}p_3(n),\quad n\geq M_c.
\end{equation}
Also, for $n\geq 58$, we have
\begin{equation}\label{b-M-3}
  0.0172p_3(n)<M_3(a,3;n)< 0.6495p_3(n).
\end{equation}
\end{prop}
\pf
Defining $R_3(a,c,n):=G_1+G_2+O(n^\varepsilon)$ and combining Lemmas \ref{lem-bound-G}, \ref{lem-bound-error} and \ref{lem-bound-p3}, for odd $c$ and positive integer $n$, we have
\begin{align}\label{R-p3-rat}
 \left|\frac{R_3(a,c,n)}{p_3(n)}\right|
 &\leq\left[\frac{1.0816c^{\frac 32}\cdot n^{\frac 12}}{8n-1}\cdot e^{\frac{\mu(n)}{c}}+\frac{\frac c2\cdot n^{\frac 12}}{8n-1}\cdot e^{2\sqrt{2\delta_0}\mu(n)}\right.
 \notag\\
 &\quad\left.+(392785.264c+178.307)n^{\frac 18}\right]\cdot\frac{(8n-1)^{\frac 98}}{2(1-\frac{3}{\mu(n)})}\cdot e^{-\mu(n)}.
\end{align}
For $n\geq2$,
\[\frac{1}{1-\frac{3}{\mu(n)}}\leq\frac{1}{1-\frac{3}{\mu(2)}}\leq1.973.\]
Plugging this into \eqref{R-p3-rat}, for odd $c$ and integer $n\geq2$, we have
\begin{align}\label{b-R-P}
 \left|\frac{R_3(a,c,n)}{p_3(n)}\right|
 &\leq 1.3838c^{\frac 32} n^{\frac 58} e^{(\frac{1}{c}-1)\mu(n)}+0.6397cn^{\frac 58} e^{(2\sqrt{2\delta_0}-1)\mu(n)}
 \notag\\
 &\quad+(4.0201\times10^6c+1824.9112)n^{\frac 54} e^{-\mu(n)}.
\end{align}
We take the sum of all three coefficients, the highest order exponential, and the highest power of $n$ from the three terms and put them
together in one term. Since $2\sqrt{2\delta_0}>\frac 1c$ for $c\geq4$, for $n\geq2$, we have
\begin{align*}
 \left|\frac{R_3(a,c,n)}{p_3(n)}\right|
 \leq(4.0201\times10^6c^{\frac 32}+1824.9112)n^{\frac 54} e^{(2\sqrt{2\delta_0}-1)\mu(n)}=:R_c(n).
\end{align*}
This together with \eqref{M3-formula} and the definition of $R_3(a,c,n)$ gives
\begin{equation}\label{M-b-R}
  \left(\frac 1c-R_c(n)\right)p_3(n)<M_3(a,c;n)<\left(\frac 1c+R_c(n)\right)p_3(n)
\end{equation}
for $c\geq4$ and $n\geq2$.

For $c\geq4$ and $n\geq M_c$, we claim that
\begin{equation*}
  R_c(n)<\frac{1}{2c},
\end{equation*}
which gives \eqref{b-M-4} directly and is equivalent to prove that
\begin{equation}\label{Rc-1}
  \frac{e^{(1-2\sqrt{2\delta_0})\mu(n)}}{n^{\frac 54}}>2c\cdot(4.0201\times10^6c^{\frac 32}+1824.9112).
\end{equation}
In addition, we recall the following inequality \cite{Mitrinovic-1964}
\begin{equation*}
  e^x>\left(1+\frac xy\right)^y,\quad x,y>0.
\end{equation*}
Hence, letting $x=(1-2\sqrt{2\delta_0})\mu(n)$ and $y=3$, we have
\begin{align}\label{Rc-2}
\frac{e^{(1-2\sqrt{2\delta_0})\mu(n)}}{n^{\frac 54}}
&> \frac{1}{n^{\frac 54}}\left(1+\frac{(1-2\sqrt{2\delta_0})\mu(n)}{3}\right)^3
\notag\\
&>\frac{\mu(n)^3(1-2\sqrt{2\delta_0})^3}{3^3n^{\frac 54}}.
\end{align}
In view of \eqref{Rc-1}, \eqref{Rc-2} and the definitions of $\mu(n)$ and $\delta_0$, it is sufficient to verify
\begin{equation}\label{Rc-3}
  \frac{(2n-\frac 14)^{\frac 32}}{n^{\frac 54}}>\frac{54c\cdot(4.0201\times10^6c^{\frac 32}+1824.9112)}{\left(\pi-2\pi\sqrt{\frac{1}{c^2}-\frac{1}{c}+\frac{1}{4}}\right)^3}.
\end{equation}
Since $1-\frac 1n\geq \frac 12$ for $n\geq2$, we have
\begin{align*}
  \frac{(2n-\frac 14)^{\frac 32}}{n^{\frac 54}}> \frac{2^{\frac 32}(n-1)^{\frac 32}}{n^{\frac 54}}=2^{\frac 32}\cdot(1-\frac 1n)^{\frac 54}\cdot(n-1)^{\frac 14}
  \geq2^{\frac 14}(n-1)^{\frac 14}.
\end{align*}
Therefore \eqref{Rc-3} holds for $n\geq M_c$.

We next consider the case $c=3$. In this case, the highest order exponential of \eqref{b-R-P} is $-\frac{2}{3}\mu(n)$. Setting $c=3$ in \eqref{b-R-P}, for $n\geq2$, we similarly have
\begin{align*}
 \left|\frac{R_3(a,3,n)}{p_3(n)}\right|
 \leq1.2063\times10^7n^{\frac 54} e^{-\frac{2}{3}\mu(n)}=:R(n).
\end{align*}
For $n\geq58$, we compute that
\begin{equation*}
  R(n)<0.3161.
\end{equation*}
Replacing $R_c(n)$ with $R(n)$ and letting $c=3$ in \eqref{M-b-R}, we arrive at \eqref{b-M-3}.
\qed

With the bounds of $M_3(a,c;n)$, we can prove Theorem \ref{ineq-sum} now.

{\noindent\it{Proof of Theorem \ref{ineq-sum}.}}
For $c\geq4$, let $n_1=\lambda n_2$ with $\lambda\geq1$. Then by \eqref{bound-p3} and \eqref{b-M-4}, we have
\begin{align*}
  M_3(a,c;n_1)M_3(a,c;n_2)>\frac{\left(1-\frac{3}{\mu(\lambda n_2)}\right)\left(1-\frac{3}{\mu(n_2)}\right)}{c^2(8\lambda n_2-1)^{\frac 98}(8n_2-1)^{\frac 98}}\cdot e^{\mu(\lambda n_2)+\mu(n_2)}
\end{align*}
and
\begin{align*}
  M_3(a,c;n_1+n_2)<\frac{\frac 97+\frac{3}{\mu(\lambda n_2+n_2)}}{c\left(8(\lambda n_2+n_2)-1\right)^{\frac 98}}\cdot e^{\mu(\lambda n_2+n_2)}
\end{align*}
for $n\geq M_c$.

Hence it suffices to show that
\begin{equation*}
  T_{n_2}(\lambda)>\log(V_{n_2}(\lambda))+\log(S_{n_2}(\lambda)),
\end{equation*}
where
\begin{align*}
  T_{n_2}(\lambda)&:=\mu(\lambda n_2)+\mu(n_2)-\mu(\lambda n_2+n_2),\\
  V_{n_2}(\lambda)&:=c\cdot\left(\frac{(8\lambda n_2-1)(8n_2-1)}{8(\lambda n_2+n_2)-1}\right)^{\frac 98},\\
  S_{n_2}(\lambda)&:=\frac{\frac 97+\frac{3}{\mu(\lambda n_2+n_2)}}{\left(1-\frac{3}{\mu(\lambda n_2)}\right)\left(1-\frac{3}{\mu(n_2)}\right)}.
\end{align*}
On the one hand, it can be shown that $T_{n_2}(\lambda)$ is increasing and $S_{n_2}(\lambda)$ is decreasing for $\lambda\geq 1$. On the other hand,
\begin{equation}\label{V-1}
  \frac{(8\lambda n_2-1)(8n_2-1)}{8(\lambda n_2+n_2)-1}\leq\frac{8\lambda n_2\cdot8n_2}{8\lambda n_2+8n_2-1}<8n_2,
\end{equation}
So it is sufficient to verify that
\begin{align}\label{4-1}
  T_{n_2}(1)&=\pi\left(2\sqrt{2n_2-\frac 14}-\sqrt{4n_2-\frac 14}\right)
  \notag\\
  &>\log\left(c(8n_2)^{\frac 98}\right)+\log\left(S_{n_2}(1)\right).
\end{align}
Since $S_{n_2}(1)$ is decreasing for $n_2\geq 2$, we have
\begin{align}\label{4-2}
\log\left(c(8n_2)^{\frac 98}\right)+\log\left(S_{n_2}(1)\right)
<\log\left(c(8n_2)^{\frac 98}\right)+\log\left(\frac{\frac 97+\frac{3}{\mu(4)}}{\left(1-\frac{3}{\mu(2)}\right)^2}\right)
<\log\left(66cn^{\frac 98}\right).
\end{align}
If $n_2\geq66c$, it is obvious that
\begin{align}\label{4-3}
\pi\left(2\sqrt{2n_2-\frac 14}-\sqrt{4n_2-\frac 14}\right)>2\log n_2^{\frac 98}=\log n_2^{\frac 98}+\log n_2^{\frac 98}>\log\left(66cn^{\frac 98}\right).
\end{align}
In view of \eqref{4-1}, \eqref{4-2} and \eqref{4-3}, we arrive at \eqref{thm-ineq} for $c\geq4$.

For $c=3$,  by \eqref{bound-p3} and \eqref{b-M-3}, repeating the above steps we have
\begin{align*}
  M_3(a,3;n_1)M_3(a,3;n_2)>(0.0344)^2\cdot\frac{\left(1-\frac{3}{\mu(\lambda n_2)}\right)\left(1-\frac{3}{\mu(n_2)}\right)}{(8\lambda n_2-1)^{\frac 98}(8n_2-1)^{\frac 98}}\cdot e^{\mu(\lambda n_2)+\mu(n_2)}
\end{align*}
and
\begin{align*}
  M_3(a,3;n_1+n_2)<\frac{1.299\times\left(\frac 97+\frac{3}{\mu(\lambda n_2+n_2)}\right)}{\left(8(\lambda n_2+n_2)-1\right)^{\frac 98}}\cdot e^{\mu(\lambda n_2+n_2)}
\end{align*}
for $n\geq 58$.

Hence it suffices to show that
\begin{equation*}
  T_{n_2}(\lambda)>\log(\overline{V}_{n_2}(\lambda))+\log(S_{n_2}(\lambda)),
\end{equation*}
where
\begin{equation*}
 \overline{V}_{n_2}(\lambda):=\frac{1.299}{(0.0344)^2}\cdot\left(\frac{(8\lambda n_2-1)(8n_2-1)}{8(\lambda n_2+n_2)-1}\right)^{\frac 98}.
\end{equation*}
Using \eqref{V-1}, we have
\begin{equation*}
  \overline{V}_{n_2}(\lambda)<11389n_2^{\frac 98}.
\end{equation*}
So it suffices to verify that
\begin{align*}
  T_{n_2}(1)=\pi\left(2\sqrt{2n_2-\frac 14}-\sqrt{4n_2-\frac 14}\right)>\log\left(11389n_2^{\frac 98}\right)+\log(S_{n_2}(1)),
\end{align*}
which holds for $n\geq 28$. It is routine to check that \eqref{thm-ineq} is true for $c=3$ and $4\leq n_1,n_2\leq57$.
\qed

\vspace{0.5cm}
 \baselineskip 15pt
{\noindent\bf\large{\ Acknowledgements}} \vspace{7pt} \par
\noindent This work was supported by the National Natural Science Foundation of China (Grant Nos. 12001182 and 12171487), the Fundamental Research Funds for the Central Universities (Grant No. 531118010411), and the Hunan Provincial Natural Science Foundation of China (Grant No. 2021JJ40037).

\end{document}